# CENTRAL LIMIT THEOREM FOR LINEAR EIGENVALUE STATISTICS OF RANDOM MATRICES WITH INDEPENDENT ENTRIES

By A. Lytova and L. Pastur

*Mathematical Division B. Verkin Institute for Low Temperature Physics,
Kharkov, Ukraine*

We consider $n \times n$ real symmetric and Hermitian Wigner random matrices $n^{-1/2}W$ with independent (modulo symmetry condition) entries and the (null) sample covariance matrices $n^{-1}X^*X$ with independent entries of $m \times n$ matrix $X$. Assuming first that the 4th cumulant (excess) $\kappa_4$ of entries of $W$ and $X$ is zero and that their 4th moments satisfy a Lindeberg type condition, we prove that linear statistics of eigenvalues of the above matrices satisfy the central limit theorem (CLT) as $n \to \infty$, $m \to \infty$, $m/n \to c \in [0, \infty)$ with the same variance as for Gaussian matrices if the test functions of statistics are smooth enough (essentially of the class $\mathbf{C}^5$). This is done by using a simple "interpolation trick" from the known results for the Gaussian matrices and the integration by parts, presented in the form of certain differentiation formulas. Then, by using a more elaborated version of the techniques, we prove the CLT in the case of nonzero excess of entries again for essentially $\mathbb{C}^5$ test function. Here the variance of statistics contains an additional term proportional to $\kappa_4$. The proofs of all limit theorems follow essentially the same scheme.

**1. Introduction.** The central limit theorem (CLT) is an important and widely used ingredient of asymptotic description of stochastic objects. In the random matrix theory, more precisely, in its part that deals with asymptotic distribution of eigenvalues $\{\lambda_l^{(n)}\}_{l=1}^n$ of random matrices of large size $n$, natural objects to study are linear eigenvalue statistics, defined via test functions $\varphi : \mathbb{R} \to \mathbb{C}$ as

$$(1.1) \qquad \mathcal{N}_n[\varphi] = \sum_{l=1}^n \varphi(\lambda_l^{(n)}).$$

---

---

The question of fluctuations of linear eigenvalue statistics of random matrices was first addressed by Arharov [3], who announced the convergence in probability of any finite collection of properly normalized traces of powers of sample covariance matrices in the case where the numbers of rows and columns of the data matrix are of the same order [see formulas (4.1), (4.2) and (4.7) below]. The result was restated and proved by Jonsson [16]. However, the explicit form of the variance of the limiting Gaussian law was not given in [3] and [16]. In 1975 Girko considered the CLT for the traces of resolvent of the Wigner and the sample covariance matrices by combining the Stieltjes transform and the martingale techniques (see [12] for results and references). In particular, an expression for the variance of the limiting Gaussian laws was given, although the expression is much less explicit than our formulas (3.92) and (4.65) for $\varphi(\lambda) = (\lambda - z)^{-1}$, $\Im z \neq 0$. In the last decade a number of results on the CLT for linear eigenvalue statistics of various classes of random matrices has been obtained (see [2, 5, 7, 8, 10, 11, 13, 15, 16, 17, 18, 28, 29, 30] and [31] and the references therein).

A rather unusual property of linear eigenvalue statistics is that their variance remains bounded as $n \to \infty$ for test functions with bounded derivative. This has to be compared with the case of linear statistics of independent and identically distributed random variables $\{\xi_l^{(n)}\}_{l=1}^n$, where the variance is linear in $n$ for any bounded test function. This fact is an important element of the ideas and techniques of the proof of the CLT for

$$(1.2) \qquad (\mathcal{N}_n[\varphi] - \mathbf{E}\{\mathcal{N}_n[\varphi]\})/(\mathbf{Var}\{\mathcal{N}_n[\varphi]\})^{1/2},$$

viewed as a result of addition of large numbers of small terms (see, e.g., [14], Chapter 18). On the other hand, since the variance of linear statistics of eigenvalues of many random matrices is bounded in $n$, the CLT, if any, has to be valid for statistics (1.1) themselves (i.e., without any $n$-dependent normalizing factor in front), resulting from a rather subtle cancelation between the terms of the sum. One can also imagine that the cancelation is not always the case, and indeed it was shown in [24] that the CLT is not necessarily valid for so-called Hermitian matrix models, for which non-Gaussian limiting laws appear in certain cases even for real analytic test functions.

In this paper we prove the CLT for linear eigenvalue statistics of two classes of random matrices: the Wigner matrices $n^{-1/2}W$, where $W$ are $n \times n$ real symmetric random matrices with independent (modulo symmetry conditions) entries (typically $n$-independent) and the matrices $n^{-1}X^TX$, where $X$ are $m \times n$ matrices with independent (and also typically $n$-independent) entries. We will refer to these matrices as the Wigner and the sample covariance matrices, respectively. The case, where the entries of $W$ are Gaussian and the probability law of $W$ is orthogonal invariant, is known as the Gaussian orthogonal ensemble (GOE). Likewise, the case, where the entries of $X$ are i.i.d. Gaussian, is known as the (null or white) Wishart ensemble. In particular, the Wishart ensemble has been used in statistics since the 30s as an important element of the sample covariance analysis, the principal component analysis first of all, in the asymptotic regime $n \to \infty, m < \infty$ (see, e.g., [21] and the references therein). The eigenvalue distribution of these matrices for $n \to \infty, m \to \infty, m/n \to c \in [0, \infty)$, that is, an analog of the law of large numbers for $n^{-1}\mathcal{N}_n[\varphi]$, was found in [20].

The CLT for certain linear eigenvalue statistics of the Wigner and the sample covariance matrices were also considered in recent papers [2, 5] and [8]. In [2] the Wigner and the sample covariance matrices (in fact, more general matrices) and linear eigenvalues statistics for polynomial test functions were studied by using a considerable amount of nontrivial combinatorics, that is, in fact, a version of the moment method of proof of the CLT. This requires the existence of all moments of entries and certain conditions on their growth as their order tends to infinity. The conditions were then relaxed for differentiable test functions under the additional assumption that the probability law of entries satisfies a concentration inequality of the Poincaré type.

Related results are obtained in [8] for a special class of Wigner matrices, whose entries have the form

$$
(1.3) \qquad W_{jk} = u(\widehat{W}_{jk}), \qquad u \in C^2(\mathbb{R}),
$$
$$
\sup_{x \in \mathbb{R}} |u'(x)| < \infty, \qquad \sup_{x \in \mathbb{R}} |u''(x)| < \infty,
$$

where $\{\widehat{W}_{jk}\}_{1 \leq j \leq k \leq n}$ are the independent standard ($\mathbf{E}\{\widehat{W}_{jk}\} = 0, \mathbf{E}\{\widehat{W}_{jk}^2\} = 1$) Gaussian random variables. For these, rather "close" to the Gaussian, random matrices a nice bound for the total variation distance between the laws of their linear eigenvalue statistics and the corresponding Gaussian random variable was given. The bound is then used to prove the CLT for linear eigenvalue statistics with entire test functions without explicit formula for the variance with a possibility of extension to $C^1$ functions and also for certain polynomials of growing with $n$ degree, as in [29].

In [5] the real symmetric and Hermitian sample covariance matrices (in fact, more general matrices) were studied, assuming that the entries $X_{\alpha j}$, $\alpha = 1, \ldots, m$, $j = 1, \ldots, n$, of $X$ are such that

$$
(1.4) \qquad \mathbf{E}\{X_{11}\} = 0, \qquad \mathbf{E}\{X_{11}^2\} = 1, \qquad \mathbf{E}\{X_{11}^4\} = 3
$$

in the real symmetric case and

$$
(1.5) \quad \mathbf{E}\{X_{11}\} = \mathbf{E}\{X_{11}^2\} = 0, \qquad \mathbf{E}\{|X_{11}|^2\} = 1, \qquad \mathbf{E}\{|X_{11}|^4\} = 2
$$

in the Hermitian case. Under these conditions the CLT for linear eigenvalue statistics with real analytic test functions was proved.

Conditions (1.4) and (1.5) mean that the fourth cumulant (known in statistics as the excess) of entries is zero. On the other hand, it was shown in [18] that for $\varphi(\lambda) = (\lambda - z)^{-1}, \Im z \neq 0$ the variance of the corresponding linear statistic (the trace of resolvent) of Wigner matrices contains the fourth cumulant of entries. Thus, even in the class of real analytic test functions one can expect more in the case of nonzero fourth cumulant of entries.

The requirement of the real analyticity of test functions results from the use of the Stieltjes transform of the eigenvalue counting measure as a basic characteristic (moment generating) function. The Stieltjes transform was introduced in the random matrix studies in [20] and since then proved to be useful in a number of problems of the field (see, e.g., [4, 12, 18] and [23] and the references therein). We found, however, that while studying the CLT of the above ensembles it is more convenient to use as a basic characteristic function not the collection of moments or the Stieltjes transform but the Fourier transform of the eigenvalue counting measure, that is, the standard characteristic function of probability theory. This allows us to prove the CLT for linear eigenvalue statistics with sufficiently regular (essentially $C^5$) test functions (but not real analytic as in [5]) and assuming the existence of the fourth moments of entries satisfying a Lindeberg type condition (but not all the moments or the Poincaré type inequality as in [2], conditions (1.3) as in [8], or conditions (1.4) and (1.5) as in [5]). Besides, all proofs follow the same scheme based on the systematic use of rather simple means: the Fourier transform, certain differential formulas, that is, a version of integration by parts [see (2.20) and (3.6)], and an "interpolation trick" [25], which allows us to relate the asymptotic properties of a number of important quantities for general entries and those for the Gaussian entries. For both classes of random matrices we prove first the CLT for matrices with Gaussian entries (the GOE and the Wishart ensemble, see Theorems 2.2 and 4.2), then consider matrices with zero excess of entries, where the CLT can be obtained practically directly from that for the GOE and the Wishart ensemble by using an "interpolation" trick (see the proof of Theorems 3.3, 3.4 and 4.3). Finally, the proofs in the general case of nonzero excess of entries essentially follow those of the GOE and the Wishart cases and use again the interpolation trick that makes the proofs shorter and simpler (see Theorems 3.6 and 4.5).

The paper is organized as follows. In Section 2 we present the basics of our approach by proving the CLT for linear eigenvalue statistics in a technically simple case of the GOE (for other proofs see, e.g., [7, 12, 13] and [15] and the references therein). In Section 3 we consider the Wigner matrices and in Section 4 the sample covariance matrices. We confine ourselves to real symmetric matrices, although our results as well as the main ingredients of the proofs remain valid in the Hermitian case with natural modifications.

Throughout the paper we write the integrals without limits for the integrals over the whole real axis.

## 2. Gaussian orthogonal ensemble.

2.1. *Generalities.* We recall first several technical facts that will be often used below. We start from the generalized Fourier transform, in fact, the $\pi/2$ rotated Laplace transform (see, e.g., [32], Sections 1.8–9).

PROPOSITION 2.1. *Let $f:\mathbb{R}_+ \to \mathbb{C}$ be locally Lipshitzian and such that for some $\delta > 0$*

$$(2.1) \qquad \sup_{t \geq 0} e^{-\delta t} |f(t)| < \infty$$

*and let $\widetilde{f}:\{z \in \mathbb{C} : \Im z < -\delta\} \to \mathbb{C}$ be its generalized Fourier transform*

$$(2.2) \qquad \widetilde{f}(z) = i^{-1} \int_0^\infty e^{-izt} f(t)\, dt.$$

*The inversion formula is given by*

$$(2.3) \qquad f(t) = \frac{i}{2\pi} \int_L e^{izt} \widetilde{f}(z)\, dz, \qquad t \geq 0,$$

*where $L = (-\infty - i\varepsilon, \infty - i\varepsilon)$, $\varepsilon > \delta$, and the principal value of the integral at infinity is used.*

*Denote for the moment the correspondence between functions and their generalized Fourier transforms as $f \leftrightarrow \widetilde{f}$. Then we have the following:*

(i) *$f'(t) \leftrightarrow i(f(+0) + z\widetilde{f}(z))$;*
(ii) *$\int_0^t f(\tau)\, d\tau \leftrightarrow (iz)^{-1} \widetilde{f}(z)$;*
(iii) *$\int_0^t f_1(t - \tau) f_2(\tau)\, d\tau := (f_1 * f_2)(t) \leftrightarrow i\widetilde{f_1}(z)\widetilde{f_2}(z)$;*
(iv) *if $P$, $Q$ and $R$ are locally Lipshitzian, satisfy (2.1), and*

$$(2.4) \qquad 1 + i\widetilde{Q}(z) \neq 0, \qquad \Im z < 0,$$

*then the equation*

$$(2.5) \qquad P(t) + \int_0^t Q(t - t_1) P(t_1)\, dt_1 = R(t), \qquad t \geq 0,$$

*has a unique locally Lipshitzian solution*

$$(2.6) \qquad P \leftrightarrow \widetilde{R}(1 + i\widetilde{Q})^{-1}$$

*or*

$$(2.7) \qquad P(t) = -i \int_0^t T(t - t_1) R(t_1)\, dt_1,$$

*where*

$$(2.8) \qquad T \leftrightarrow (1 + i\widetilde{Q})^{-1}.$$

*In particular, if $R(t)$ is differentiable, $R(0) = 0$, and*

$$(2.9) \qquad\qquad Q(t) = \int_0^t Q_1(s) \, ds,$$

*then the equation*

$$(2.10) \qquad P(t) + \int_0^t dt_1 \int_0^{t_1} Q_1(t_1 - t_2) P(t_2) \, dt_2 = R(t), \qquad t \geq 0,$$

*has a unique locally Lipshitzian solution*

$$(2.11) \qquad\qquad P(t) = - \int_0^t T_1(t - t_1) R'(t_1) \, dt_1,$$

*where*

$$(2.12) \qquad\qquad T_1 \leftrightarrow (z + \widetilde{Q}_1)^{-1}$$

*provided by*

$$(2.13) \qquad\qquad z + \widetilde{Q}_1(z) \neq 0, \qquad \Im z < 0.$$

The next proposition presents a simple fact of linear algebra:

PROPOSITION 2.2 (Duhamel formula).    *Let $M_1$, $M_2$ be $n \times n$ matrices and $t \in \mathbb{R}$. Then we have*

$$(2.14) \qquad e^{(M_1 + M_2)t} = e^{M_1 t} + \int_0^t e^{M_1(t-s)} M_2 e^{(M_1 + M_2)s} \, ds.$$

Consider now a real symmetric matrix $M = \{M_{jk}\}_{j,k=1}^n$ and set

$$(2.15) \qquad\qquad U(t) = e^{itM}, \qquad t \in \mathbb{R}.$$

Then $U(t)$ is a symmetric unitary matrix, possessing the properties

$$(2.16) \qquad U(t)U(s) = U(t+s),$$

$$\|U(t)\| = 1, \qquad |U_{jk}(t)| \leq 1, \qquad \sum_{j=1}^n |U_{jk}(t)|^2 = 1.$$

The Duhamel formula allows us to obtain the derivatives of $U(t)$ with respect to the entries $M_{jk}$, $j, k = 1, \dots, n$, of $M$:

$$(2.17) \quad D_{jk} U_{ab}(t) = i\beta_{jk}[(U_{aj} * U_{bk})(t) + (U_{bj} * U_{ak})(t)], \qquad D_{jk} = \partial/\partial M_{jk},$$

*where*

$$(2.18) \qquad \beta_{jk} = (1 + \delta_{jk})^{-1} = \begin{cases} 1/2, & j = k, \\ 1, & j \neq k, \end{cases}$$

and the symbol "$*$" is defined in Proposition 2.1(iii). Iterating (2.17) and using (2.16), we obtain the bound

$$(2.19) \qquad |D_{jk}^l U_{ab}(t)| \le c_l |t|^l,$$

where $c_l$ is an absolute constant for every $l$.

The next proposition presents certain facts on Gaussian random variables.

PROPOSITION 2.3. *Let $\xi = \{\xi_l\}_{l=1}^p$ be independent Gaussian random variables of zero mean, and $\Phi : \mathbb{R}^p \to \mathbb{C}$ be a differentiable function with polynomially bounded partial derivatives $\Phi_l', l = 1, \ldots, p$. Then we have*

$$(2.20) \qquad \mathbf{E}\{\xi_l \Phi(\xi)\} = \mathbf{E}\{\xi_l^2\}\mathbf{E}\{\Phi_l'(\xi)\}, \qquad l = 1, \ldots, p,$$

*and*

$$(2.21) \qquad \mathbf{Var}\{\Phi(\xi)\} \le \sum_{l=1}^p \mathbf{E}\{\xi_l^2\}\mathbf{E}\{|\Phi_l'(\xi)|^2\}.$$

The first formula is a version of the integration by parts. The second is a version of the Poincaré inequality (see, e.g., [6]).

Next is the definition of the Gaussian orthogonal ensemble. This is a real symmetric $n \times n$ random matrix

$$(2.22) \qquad M = n^{-1/2}W, \qquad W = \{W_{jk} \in \mathbb{R}, W_{jk} = W_{kj}\}_{j,k=1}^n,$$

defined by the probability law

$$(2.23) \qquad Z_{n1}^{-1} e^{-\operatorname{Tr} W^2 / 4w^2} \prod_{1 \le j \le k \le n} dW_{jk},$$

where $Z_{n1}$ is the normalization constant. Since

$$\operatorname{Tr} W^2 = \sum_{1 \le j \le n} W_{jj}^2 + 2 \sum_{1 \le j < k \le n} W_{jk}^2,$$

the above implies that $\{W_{jk}\}_{1 \le j \le k \le n}$ are independent Gaussian random variables such that

$$(2.24) \qquad \mathbf{E}\{W_{jk}\} = 0, \qquad \mathbf{E}\{W_{jk}^2\} = w^2(1 + \delta_{jk}).$$

Here is a useful bound for linear eigenvalue statistics of the GOE [9, 23]:

PROPOSITION 2.4. *Let $M$ be the GOE matrix (2.22)–(2.24) and $\mathcal{N}_n[\varphi]$ be its linear eigenvalue statistic (1.1), corresponding to a differentiable test function. Then*

$$(2.25) \qquad \mathbf{Var}\{\mathcal{N}_n[\varphi]\} \le 2w^2 \mathbf{E}\{n^{-1} \operatorname{Tr} \varphi'(M)(\varphi'(M)^*)\}$$

$$(2.26) \qquad \le 2w^2 \Big(\sup_{\lambda \in \mathbb{R}} |\varphi'(\lambda)|\Big)^2.$$

PROOF. The spectral theorem for real symmetric matrices and (1.1) imply that

$$\mathcal{N}_n[\varphi] = \operatorname{Tr} \varphi(M). \tag{2.27}$$

Thus, we can apply (2.21) to $\Phi(M) = \operatorname{Tr} \varphi(M)$, viewing it as a differentiable function of the independent Gaussian random variables $M_{jk} = n^{-1/2} W_{jk}$, $1 \le j \le k \le n$, satisfying (2.24). By using (2.21), (2.24) and the formula [see (2.17)]

$$D_{jk} \operatorname{Tr} \varphi(M) = 2\beta_{jk} \varphi'_{jk}(M), \tag{2.28}$$

we obtain

$$\mathbf{Var}\{\mathcal{N}_n[\varphi]\} \le w^2 n^{-1} \sum_{1 \le j \le k \le n} 4(1 + \delta_{jk}) \beta_{jk}^2 \mathbf{E}\{|\varphi'_{jk}(M)|^2\}$$

$$= 2w^2 n^{-1} \sum_{j,k=1}^{n} \mathbf{E}\{|\varphi'_{jk}(M)|^2\} = 2w^2 \mathbf{E}\{n^{-1} \operatorname{Tr} \varphi'(M)(\varphi'(M))^*\}.$$

This yields (2.25). Using it and the inequalities

$$|\operatorname{Tr} A| \le n\|A\|, \qquad \|\psi(B)\| \le \sup_{\lambda \in \mathbb{R}} |\psi(\lambda)|, \tag{2.29}$$

valid for any normal matrix $A$, Hermitian matrix $B$, and $\psi : \mathbb{R} \to \mathbb{C}$, we obtain (2.26). $\square$

We recall now an analog of the law of large numbers for linear eigenvalue statistics of the GOE (see, e.g., [4, 12] and [23] and the references therein).

THEOREM 2.1. *Let $M$ be the GOE matrix (2.22)–(2.24), and $\mathcal{N}_n[\varphi]$ be a linear statistic of its eigenvalues (1.1). Then we have for any bounded and continuous $\varphi : \mathbb{R} \to \mathbb{C}$ with probability 1*

$$\lim_{n \to \infty} n^{-1} \mathcal{N}_n[\varphi] = \int \varphi(\lambda) N_{\mathrm{scl}}(d\lambda), \tag{2.30}$$

*where the measure*

$$N_{\mathrm{scl}}(d\lambda) = \rho_{\mathrm{scl}}(\lambda) \, d\lambda, \qquad \rho_{\mathrm{scl}}(\lambda) = (2\pi w^2)^{-1} (4w^2 - \lambda^2)_+^{1/2} \tag{2.31}$$

*is known as the Wigner or the semicircle law, and $x_+ = \max\{0, x\}$.*

We need below a weak version of the theorem in which the convergence with probability 1 is replaced by the convergence in mean, that is,

$$\lim_{n \to \infty} \mathbf{E}\{n^{-1} \mathcal{N}_n[\varphi]\} = \int \varphi(\lambda) N_{\mathrm{scl}}(d\lambda) \tag{2.32}$$

for any continuous and bounded $\varphi$. We outline the proof of this relation to introduce several elements of techniques used in the paper (see [23] for details).

Introduce the normalized counting (empirical) measure of eigenvalues as

$$(2.33) \qquad N_n(\Delta) = \sharp\{\lambda_l^{(n)} \in \Delta : l = 1, \ldots, n\}/n.$$

Then we have

$$\mathbf{E}\{n^{-1}\mathcal{N}_n[\varphi]\} = \int \varphi(\lambda)\mathbf{E}\{N_n(d\lambda)\},$$

hence, (2.32) is equivalent to the weak convergence of $\mathbf{E}\{N_n\}$ to $N_{\mathrm{scl}}$. Moreover, since by (2.24)

$$(2.34) \qquad \int \lambda^2 \mathbf{E}\{N_n(d\lambda)\} = \mathbf{E}\{n^{-2}\operatorname{Tr} W^2\} = (1 + n^{-1})w^2,$$

the sequence $\{\mathbf{E}\{N_n\}\}$ is tight, and it suffices to prove the vague convergence of the sequence, for instance, the convergence of the Stieltjes transforms

$$(2.35) \qquad f_n(z) = \int \frac{\mathbf{E}\{N_n(d\lambda)\}}{\lambda - z}$$

of $\mathbf{E}\{N_n\}$ for any $\Im z \neq 0$ (see, e.g., [1], Section 59) to the Stieltjes transform

$$(2.36) \qquad f(z) = (\sqrt{z^2 - 4w^2} - z)/2w^2$$

of $N_{\mathrm{scl}}$, where $\sqrt{z^2 - 4w^2}$ is defined by the asymptotic $\sqrt{z^2 - 4w^2} = z + O(z^{-1})$, $z \to \infty$.

It follows from the definition of $N_n$ that

$$(2.37) \qquad f_n(z) = \mathbf{E}\{n^{-1}\operatorname{Tr} G(z)\},$$

where

$$G(z) = (M - z)^{-1}, \qquad \Im z \neq 0,$$

is the resolvent of $M$. We will need the resolvent identity

$$(2.38) \quad (A + B - z)^{-1} - (A - z)^{-1} = -(A + B - z)^{-1}B(A - z)^{-1},$$

its implication

$$(2.39) \qquad \frac{d}{d\varepsilon}(A + \varepsilon B - z)^{-1}\Big|_{\varepsilon=0} = -(A - z)^{-1}B(A - z)^{-1}$$

and the bounds

$$(2.40) \qquad \|(A - z)^{-1}\| \leq |\Im z|^{-1}, \qquad |((A - z)^{-1})_{jk}| \leq |\Im z|^{-1},$$

valid for real symmetric matrices $A$ and $B$.

It follows from (2.38) for $A = 0, B = M$ that

$$f_n(z) = -\frac{1}{z} + \frac{2}{zn^{3/2}}\mathbf{E}\left\{\sum_{1 \leq j \leq k \leq n} \beta_{jk}W_{jk}G_{jk}(z)\right\},$$

where $\beta_{jk}$ is defined in (2.18). Since $W_{jk}, 1 \leq j \leq k \leq n$, are independent Gaussian variables, we can write, in view of (2.20) and (2.24),

$$(2.41) \qquad f_n(z) = -\frac{1}{z} + \frac{2w^2}{zn^2}\mathbf{E}\left\{\sum_{1 \leq j \leq k \leq n} D_{jk}G_{jk}(z)\right\},$$

where $D_{jk}$ is defined in (2.17). It follows from (2.39) that [cf. (2.17)]

$$(2.42) \qquad D_{jk}G_{ab}(z) = -\beta_{jk}(G_{aj}(z)G_{kb}(z) + G_{ak}(z)G_{jb}(z)).$$

This allows us to write (2.41) as

$$(2.43) \quad f_n(z) = -z^{-1} - w^2 z^{-1}\mathbf{E}\{g_n^2(z)\} - w^2 z^{-1}\mathbf{E}\{n^{-2}\operatorname{Tr}G^2(z)\},$$

where

$$(2.44) \qquad\qquad g_n(z) = n^{-1}\operatorname{Tr}G(z).$$

By using (2.26) with $\varphi(\lambda) = (\lambda - z)^{-1}$, we find that

$$(2.45) \qquad\qquad \mathbf{Var}\{g_n(z)\} \leq \frac{2w^2}{n^2|\Im z|^4},$$

hence,

$$|\mathbf{E}\{g_n^2(z)\} - f_n^2(z)| \leq \mathbf{Var}\{g_n(z)\} \leq \frac{2w^2}{n^2|\Im z|^4}.$$

Besides, (2.29) and (2.40) imply that $|\mathbf{E}\{n^{-2}\operatorname{Tr}G^2(z)\}| \leq 1/n|\Im z|^2$, and (2.35) implies that $|f_n| \leq 1/|\Im z|$.

In view of the above bounds, the sequence $\{f_n\}$ is compact with respect to the uniform convergence on any compact set $K \subset \mathbb{C} \setminus \mathbb{R}$ and the uniform on $K$ limit $f$ of any convergent subsequence of $\{f_n\}$ satisfies the quadratic equation

$$(2.46) \qquad\qquad f(z) = -z^{-1} - w^2 z^{-1}f^2(z),$$

following from (2.43). In addition, we have, by (2.35), $\Im f_n(z)\Im z > 0$, thus, $\Im f(z)\Im z \geq 0$. Now it is elementary to check that the unique solution of the above quadratic equation that satisfies this condition is (2.36).

2.2. *Central limit theorem for linear eigenvalue statistics of differentiable test functions.* The CLT for the GOE was proved in a number of works (see [7, 13] and [23] and the references therein). We present below a proof, whose strategy dates back to [18] and is used in the proofs of other CLT of the paper.

THEOREM 2.2. *Let $\varphi : \mathbb{R} \to \mathbb{R}$ be a bounded function with bounded derivative, and $\mathcal{N}_n[\varphi]$ be the corresponding linear eigenvalue statistic (1.1) of the GOE (2.22)–(2.24). Then the random variable*

$$(2.47) \qquad \mathcal{N}_n^\circ[\varphi] = \mathcal{N}_n[\varphi] - \mathbf{E}\{\mathcal{N}_n[\varphi]\}$$

*converges in distribution to the Gaussian random variable with zero mean and variance*

$$(2.48) \quad V_{\mathrm{GOE}}[\varphi] = \frac{1}{2\pi^2} \int_{-2w}^{2w} \int_{-2w}^{2w} \left( \frac{\Delta \varphi}{\Delta \lambda} \right)^2 \frac{4w^2 - \lambda_1 \lambda_2}{\sqrt{4w^2 - \lambda_1^2}\sqrt{4w^2 - \lambda_2^2}} \, d\lambda_1 \, d\lambda_2,$$

*where*

$$(2.49) \qquad \Delta \varphi = \varphi(\lambda_1) - \varphi(\lambda_2), \qquad \Delta \lambda = \lambda_1 - \lambda_2.$$

PROOF. By the continuity theorem for characteristic functions, it suffices to show that if

$$(2.50) \qquad Z_n(x) = \mathbf{E}\{e^{ix\mathcal{N}_n^\circ[\varphi]}\},$$

then for any $x \in \mathbb{R}$

$$(2.51) \qquad \lim_{n \to \infty} Z_n(x) = Z(x),$$

where

$$(2.52) \qquad Z(x) = \exp\{-x^2 V_{\mathrm{GOE}}[\varphi]/2\}.$$

We obtain first (2.52), hence the theorem, for a certain class of test functions and then we extend the theorem to the bounded functions with bounded derivative, by using a standard approximation procedure [see also Remark 2.1(2) for a wider class].

Assume then that $\varphi$ admits the Fourier transform

$$(2.53) \qquad \widehat{\varphi}(t) = \frac{1}{2\pi} \int e^{-it\lambda} \varphi(\lambda) \, d\lambda,$$

satisfying the condition

$$(2.54) \qquad \int (1 + |t|^2) |\widehat{\varphi}(t)| \, dt < \infty.$$

Following the idea of [18], we obtain (2.52) by deriving the equation

$$(2.55) \qquad Z(x) = 1 - V_{\mathrm{GOE}}[\varphi] \int_0^x y Z(y) \, dy.$$

The equation is uniquely soluble in the class of bounded continuous functions and its solution is evidently (2.52).

It follows from (2.50) that

$$(2.56) \qquad Z_n'(x) = i\mathbf{E}\{\mathcal{N}_n^\circ[\varphi] e^{ix\mathcal{N}_n^\circ[\varphi]}\}.$$

This, the Schwarz inequality and (2.26) yield

$$|Z_n'(x)| \leq \sqrt{2} w \Big( \sup_{\lambda \in \mathbb{R}} |\varphi'(\lambda)| \Big).$$

Since $Z_n(0) = 1$, we have the equality

$$(2.57) \qquad Z_n(x) = 1 + \int_0^x Z_n'(y) \, dy,$$

showing that it suffices to prove that any converging subsequences $\{Z_{n_j}\}$ and $\{Z_{n_j}'\}$ satisfy

$$(2.58) \qquad \lim_{n_j \to \infty} Z_{n_j}(x) = Z(x), \qquad \lim_{n_j \to \infty} Z_{n_j}'(x) = -x V_{\mathrm{GOE}} Z(x).$$

Indeed, if yes, then (2.58), (2.57) and the dominated convergence theorem imply (2.55), hence (2.52).

The Fourier inversion formula

$$(2.59) \qquad \varphi(\lambda) = \int e^{it\lambda} \widehat{\varphi}(t) \, dt$$

and (2.27) yield for (2.47)

$$(2.60) \qquad \mathcal{N}_n^\circ[\varphi] = \int \widehat{\varphi}(t) u_n^\circ(t) \, dt,$$

where

$$(2.61) \qquad u_n(t) = \mathrm{Tr}\, U(t), \qquad u_n^\circ(t) = u_n(t) - \mathbf{E}\{u_n(t)\},$$

and $U(t)$ defined by (2.15) with the GOE matrix $M$. Thus, we have by (2.56)

$$(2.62) \qquad Z_n'(x) = i \int \widehat{\varphi}(t) Y_n(x, t) \, dt,$$

where

$$(2.63) \qquad Y_n(x, t) = \mathbf{E}\{u_n^\circ(t) e_n(x)\}, \qquad e_n(x) = e^{ix\mathcal{N}_n^\circ[\varphi]}.$$

Since

$$(2.64) \qquad \overline{Y_n(x, t)} = Y_n(-x, -t),$$

we can confine ourselves to the half-plane $\{t \geq 0, x \in \mathbb{R}\}$.

We will prove that the sequence $\{Y_n\}$ is bounded and equicontinuous on any finite set of $\{t \geq 0, x \in \mathbb{R}\}$, and that its every uniformly converging on the set subsequence has the same limit $Y$, leading to (2.58). This proves the assertion of the theorem under condition (2.54). Indeed, let $\{Z_{n_l}\}_{l \geq 1}$ be subsequence converging to $\widetilde{Z} \neq Z$. Consider the corresponding subsequence $\{Y_{n_l}\}_{l \geq 1}$. It contains a uniformly converging sub-subsequence, whose limit is again $Y$, and this forces the corresponding subsequence of $\{Z_{n_l}\}_{l \geq 1}$ to converge to $Z$, a contradiction.

It follows from (2.25) and (2.26) with $\varphi(\lambda) = e^{it\lambda}$ and $\varphi(\lambda) = i\lambda e^{it\lambda}$ and from (2.24) that

$$\tag{2.65} \mathbf{Var}\{u_n(t)\} \leq 2w^2 t^2$$

and

$$\tag{2.66} \mathbf{Var}\{u_n'(t)\} \leq 2w^2 n^{-1} \mathbf{E}\{\mathrm{Tr}(1 + t^2 M^2)\} \leq 2w^2(1 + 2w^2 t^2).$$

This, (2.63), the Schwarz inequality, and $|e_n(x)| \leq 1$ yield

$$\tag{2.67} |Y_n(x,t)| \leq \mathbf{E}\{|u_n^\circ(t)|\} \leq \sqrt{2}w|t|$$

and

$$\tag{2.68} \left|\frac{\partial}{\partial t} Y_n(x,t)\right| \leq \mathbf{Var}^{1/2}\{u_n'(t)\} \leq \sqrt{2}w(1 + 2w^2 t^2)^{1/2}$$

and according to the Schwarz inequality, (2.26) and (2.65),

$$\tag{2.69} \begin{aligned} \left|\frac{\partial}{\partial x} Y_n(x,t)\right| &= |\mathbf{E}\{u_n^\circ(t)\mathcal{N}_n^\circ[\varphi]e_n(x)\}| \\ &\leq \mathbf{Var}^{1/2}\{u_n(t)\}\,\mathbf{Var}^{1/2}\{\mathcal{N}_n[\varphi]\} \leq 2w^2|t|\sup_{\lambda \in \mathbb{R}}|\varphi'(\lambda)|. \end{aligned}$$

We conclude that the sequence $\{Y_n\}$ is bounded and equicontinuous on any finite set of $\mathbb{R}^2$. We will prove now that any uniformly converging subsequence of $\{Y_n\}$ has the same limit $Y$, leading to (2.55), hence to (2.51) and (2.52).

It follows from the Duhamel formula (2.14) with $M_1 = 0$ and $M_2 = iM$ and (2.61) that

$$u_n(t) = n + i \int_0^t \sum_{j,k=1}^n M_{jk} U_{jk}(t_1)\, dt_1,$$

hence,

$$\tag{2.70} Y_n(x,t) = \frac{i}{\sqrt{n}} \int_0^t \sum_{j,k=1}^n \mathbf{E}\{W_{jk} U_{jk}(t_1) e_n^\circ(x)\}\, dt_1,$$

where $e_n^\circ = e_n - \mathbf{E}\{e_n\}$ or, after applying (2.20) and (2.24),

$$(2.71) \qquad Y_n(x,t) = \frac{iw^2}{n} \int_0^t \sum_{j,k=1}^n (1+\delta_{jk}) \mathbf{E}\{D_{jk}(U_{jk}(t_1)e_n^\circ(x))\} \, dt_1.$$

Now, by using (2.17) and (2.28), we obtain that [cf. (2.42)]

$$(2.72) \quad D_{jk}e_n(x) = 2i\beta_{jk}xe_n(x)\varphi'_{jk}(M) = -2\beta_{jk}xe_n(x)\int tU_{jk}(t)\widehat{\varphi}(t) \, dt,$$

where the last equality follows from [see (2.59)]

$$(2.73) \qquad \varphi'(M) = i\int \widehat{\varphi}(t)tU(t) \, dt.$$

This and (2.71) yield [cf. (2.43)]

$$Y_n(x,t) = -w^2 n^{-1} \int_0^t t_1 \mathbf{E}\{u_n(t_1)e_n^\circ(x)\} \, dt_1$$

$$- w^2 n^{-1} \int_0^t dt_1 \int_0^{t_1} \mathbf{E}\{u_n(t_2-t_1)u_n(t_2)e_n^\circ(x)\} \, dt_2$$

$$- 2w^2 x \int_0^t \mathbf{E}\{e_n(x)n^{-1}\operatorname{Tr} U(t_1)\varphi'(M)\} \, dt_1.$$

Writing

$$(2.74) \qquad \overline{v}_n(t) = n^{-1}\mathbf{E}\{u_n(t)\}$$

and

$$(2.75) \qquad u_n(t) = u_n^\circ(t) + n\overline{v}_n(t), \qquad e_n(x) = e_n^\circ(x) + Z_n(x),$$

we present the above relation for $Y_n$ as

$$(2.76) \qquad \begin{aligned} Y_n(x,t) + 2w^2 \int_0^t dt_1 \int_0^{t_1} \overline{v}_n(t_1-t_2)Y_n(x,t_2) \, dt_2 \\ = xZ_n(x)A_n(t) + r_n(x,t) \end{aligned}$$

with

$$(2.77) \qquad A_n(t) = -2w^2 \int_0^t \mathbf{E}\{n^{-1}\operatorname{Tr} U(t_1)\varphi'(M)\} \, dt_1$$

and

$$(2.78) \qquad \begin{aligned} r_n(x,t) = &-w^2 n^{-1} \int_0^t t_1 Y_n(x,t_1) \, dt_1 \\ &- w^2 n^{-1} \int_0^t dt_1 \int_0^{t_1} \mathbf{E}\{u_n^\circ(t_1-t_2)u_n^\circ(t_2)e_n^\circ(x)\} \, dt_2 \\ &- 2iw^2 xn^{-1} \int_0^t dt_1 \int t_2\widehat{\varphi}(t_2)\mathbf{E}\{u_n(t_1+t_2)e_n^\circ(x)\} \, dt_2, \end{aligned}$$

where we used (2.16) and (2.73).

It follows from the inequality $|e_n^\circ(x)| \le 2$, the Schwarz inequality, (2.54) and (2.65) that the limit

$$(2.79) \qquad \lim_{n\to\infty} r_n(x,t) = 0$$

holds uniformly on any compact of $\{t \ge 0, x \in \mathbb{R}\}$. Besides, by Theorem 2.1 the sequences $\{\overline{v}_n\}$ of (2.74) and $\{A_n\}$ of (2.77) converge uniformly on any finite interval of $\mathbb{R}$ to

$$(2.80) \qquad v(t) = \frac{1}{2\pi w^2} \int_{-2w}^{2w} e^{it\lambda} \sqrt{4w^2 - \lambda^2}\, d\lambda$$

and

$$(2.81) \qquad A(t) = -\frac{1}{\pi} \int_0^t dt_1 \int_{-2w}^{2w} e^{it_1\lambda} \varphi'(\lambda) \sqrt{4w^2 - \lambda^2}\, d\lambda.$$

The above allows us to pass to the limit $n_l \to \infty$ in (2.76) and obtain that the limit $Y$ of every uniformly converging subsequence $\{Y_{n_l}\}_{l\ge1}$ satisfies the equation

$$(2.82) \quad Y(x,t) + 2w^2 \int_0^t dt_1 \int_0^{t_1} v(t_1 - t_2) Y(x, t_2)\, dt_2 = xZ(x)A(t).$$

The equation is of the form (2.10), corresponding to $\delta = 0$ in (2.1), thus, we can use formulas (2.12) and (2.13) to write its solution.

It follows from (2.2) and (2.80) (or from the spectral theorem and Theorem 2.1) that

$$(2.83) \qquad \widehat{v} = f,$$

where $f$ is the Stieltjes transform (2.36) of the semicircle law. Thus, in our case (2.4) takes form

$$(2.84) \qquad z + 2w^2 f(z) = \sqrt{z^2 - 4w^2} \ne 0, \qquad \Im z \ne 0,$$

and

$$T_1(t) = \frac{i}{2\pi} \int_L \frac{e^{izt}}{z + 2w^2 f(z)}\, dz.$$

Replacing here the integral over $L$ by the integral over the cut $[-2w, 2w]$ and taking into account that $\sqrt{z^2 - 4w^2}$ is $\pm i\sqrt{4w^2 - \lambda^2}$ on the upper and lower edges of the cut, we find that

$$(2.85) \qquad T_1(t) = -\frac{1}{\pi} \int_{-2w}^{2w} \frac{e^{i\lambda t}\, d\lambda}{\sqrt{4w^2 - \lambda^2}}.$$

Besides, in our case the r.h.s. of (2.10) is $xZ(z)A(t)$, hence, we have by (2.11)

$$
(2.86) \quad
\begin{aligned}
Y(x,t) &= \frac{ixZ(x)}{\pi^2} \int_{-2w}^{2w} \frac{d\mu}{\sqrt{4w^2 - \mu^2}} \\
&\quad \times \int_{-2w}^{2w} \frac{e^{it\lambda} - e^{it\mu}}{(\lambda - \mu)} \varphi'(\lambda) \sqrt{4w^2 - \lambda^2} \, d\lambda, \qquad t \geq 0.
\end{aligned}
$$

According to (2.64), the same limiting expression is valid for $t \leq 0$ and $x \in \mathbb{R}$, thus, we have, in view of (2.54) and (2.62),

$$
\begin{aligned}
\lim_{n_l \to \infty} Z'_{n_l}(x) &= -\frac{xZ(x)}{\pi^2} \int_{-2w}^{2w} \frac{d\mu}{\sqrt{4w^2 - \mu^2}} \\
&\quad \times \int_{-2w}^{2w} \frac{\varphi(\lambda) - \varphi(\mu)}{(\lambda - \mu)} \varphi'(\lambda) \sqrt{4w^2 - \lambda^2} \, d\lambda.
\end{aligned}
$$

Writing

$$
\varphi'(\lambda)(\varphi(\lambda) - \varphi(\mu)) = \frac{1}{2} \frac{\partial}{\partial \lambda}(\varphi(\lambda) - \varphi(\mu))^2
$$

and integrating by parts with respect to $\lambda$, we obtain

$$
(2.87) \quad \lim_{n_l \to \infty} Z'_{n_l}(x) = -\frac{xZ(x)}{2\pi^2} \int_{-2w}^{2w} \int_{-2w}^{2w} \left(\frac{\Delta \varphi}{\Delta \lambda}\right)^2 \frac{(4w^2 - \lambda\mu) \, d\lambda \, d\mu}{\sqrt{4w^2 - \lambda^2}\sqrt{4w^2 - \mu^2}},
$$

hence (2.58), and then (2.55), thus the assertion of the theorem under the condition (2.54).

The case of bounded test functions with bounded derivative can be obtained via a standard approximation procedure. Indeed, for any bounded $\varphi$ with bounded derivative there exists a sequence $\{\varphi_k\}$ of functions, satisfying (2.54) and such that

$$
(2.88) \quad
\begin{aligned}
&\sup_{\lambda \in \mathbb{R}} |\varphi'_k(\lambda)| \leq \sup_{\lambda \in \mathbb{R}} |\varphi'(\lambda)|, \\
&\lim_{k \to \infty} \sup_{|\lambda| \leq A} |\varphi'(\lambda) - \varphi'_k(\lambda)| = 0 \qquad \forall A > 0.
\end{aligned}
$$

By the above proof we have the central limit theorem for every $\varphi_k$. Denote for the moment the characteristic functions of (2.50) and (2.52) by $Z_n[\varphi]$ and $Z[\varphi]$, to make explicit their dependence on the test function. We have then for any bounded test function with bounded derivative

$$
(2.89) \quad
\begin{aligned}
|Z_n[\varphi] - Z[\varphi]| &\leq |Z_n[\varphi] - Z_n[\varphi_k]| + |Z_n[\varphi_k] - Z[\varphi_k]| \\
&\quad + |Z[\varphi_k] - Z[\varphi]|,
\end{aligned}
$$

where the second term of the r.h.s. vanishes after the limit $n \to \infty$ in view of the above proof. It follows from (2.25), (2.33) and (2.50) that

$$|Z_n[\varphi] - Z_n[\varphi_k]| \leq |x| \mathbf{Var}^{1/2}\{\mathcal{N}_n[\varphi - \varphi_k]\}$$
$$\leq \sqrt{2}w|x|\left(\int |\varphi'(\lambda) - \varphi_k'(\lambda)|^2 \mathbf{E}\{N_n(d\lambda)\}\right)^{1/2}.$$

Now (2.88) implies that the integral on the r.h.s. is bounded by

$$2\sup_{\lambda \in \mathbb{R}} |\varphi'(\lambda)|^2 \mathbf{E}\{N_n(\mathbb{R} \setminus [-A, A])\} + \sup_{|\lambda| \leq A} |\varphi'(\lambda) - \varphi_k'(\lambda)|^2, \qquad A > 2w,$$

where the first term vanishes as $n \to \infty$ by (2.32), and the second term vanishes as $k \to \infty$ by (2.88). Besides, according to (2.52),

$$|Z[\varphi] - Z[\varphi_k]| \leq x^2 |V_{\mathrm{GOE}}[\varphi] - V_{\mathrm{GOE}}[\varphi_k]|/2$$

and taking into account the continuity of $V_{\mathrm{GOE}}$ of (2.48) with respect to the $C^1$ convergence on any interval $|\lambda| \leq A$, $A > 2w$, we find that the third term of (2.89) vanishes after the limit $k \to \infty$. Thus, we have proved the central limit theorem for bounded test functions with bounded derivative. For wider classes of test functions see [7] and [15] and Remark 2.1. □

REMARK 2.1. (1) We note that the proof of Theorem 2.2 can be easily modified to prove an analogous assertion for the Gaussian unitary ensemble of Hermitian matrices, defined by (2.22) with $W_{jk} = \overline{W_{kj}}$ and the probability law

$$Z_{n2}^{-1} e^{-\mathrm{Tr}\, W^2/2w^2} \prod_{j=1}^n dW_{jj} \prod_{1 \leq j < k \leq n}^n d\Re W_{jk}\, d\Im W_{jk}.$$

The result is given by Theorem 2.2 in which $V_{\mathrm{GOE}}$ is replaced by $V_{\mathrm{GUE}} = V_{\mathrm{GOE}}/2$.

(2) It follows from the representation of the density $\rho_n$ of $\mathbf{E}\{N_n\}$ via the Hermite polynomials (see [19], Chapters 6 and 7) or from Theorem 2.3 of [26] that

$$\rho_n(\lambda) \leq C e^{-cn\lambda^2}$$

for finite $c > 0$, $C < \infty$, and a sufficiently big $|\lambda|$. This bound and the approximation procedure of the end of proof of Theorem 2.2 allows one to extend the theorem to $C^1$ test functions whose derivative grows as $C_1 e^{c_1\lambda^2}$ for any $c_1 > 0$ and $C_1 < \infty$.

**3. Wigner ensembles.** In this section we prove the central limit theorem for linear eigenvalue statistics of the Wigner random matrices. We start from the analog of the law of large numbers, by proving that the normalized counting measure of eigenvalues converges in mean to the semicircle law. The fact has been well known since the early seventies (see [4, 12] and [22] and the references therein). We give a new proof, valid under rather general conditions of the Lindeberg type and based on a certain "interpolation" trick that is systematically used in what follows. We then pass to the CLT, proving it first for the Wigner ensembles, assuming the existence of the fourth moment of entries, their zero excess and the integrability of $(1 + |t|^5)\widehat{\varphi}$ (Theorems 3.3 and 3.4), and then proving it in the general case of an arbitrary excess (Theorem 3.6), assuming the existence of the fourth moments satisfying a Lindeberg type condition (3.57) and again the integrability of $(1 + |t|^5)\widehat{\varphi}$.

3.1. *Generalities.* We present here the definition of the Wigner ensembles and technical means that we are going to use in addition to those given in the previous section.

Wigner ensembles for real symmetric matrices can be defined as follows:

$$(3.1) \qquad M = n^{-1/2}W, \qquad W = \{W_{jk}^{(n)} \in \mathbb{R}, W_{jk}^{(n)} = W_{kj}^{(n)}\}_{j,k=1}^n$$

[cf. (2.22)], where the random variables $W_{jk}^{(n)}$, $1 \le j \le k \le n$, are independent, and

$$(3.2) \qquad \mathbf{E}\{W_{jk}^{(n)}\} = 0, \qquad \mathbf{E}\{(W_{jk}^{(n)})^2\} = (1 + \delta_{jk})w^2,$$

that is, the two first moments of the entries coincide with those of the GOE [see (2.22)–(2.24)]. In other words, the probability law of the matrix $W$ is

$$(3.3) \qquad \mathbf{P}(dW) = \prod_{1 \le j \le k \le n} F_{jk}^{(n)}(dW_{jk}),$$

where for any $1 \le j \le k \le n$, $F_{jk}^{(n)}$ is a probability measure on the real line, satisfying conditions (3.2). We do not assume in general that $F_{jk}^{(n)}$ is $n$-independent, and that $F_{jk}^{(n)}$ are the same for $1 \le j < k \le n$ and for $j = k = 1, \ldots, n$, that is, for off-diagonal and diagonal entries as in the GOE case.

Since we are going to use the scheme of the proof of Theorem 2.2 (the CLT for the GOE), we need an analog of the differentiation formula (2.20) for non-Gaussian random variables. To this end, we recall first that if a random variable $\xi$ has a finite $p$th absolute moment, $p \ge 1$, then we have the expansions

$$f(t) := \mathbf{E}\{e^{it\xi}\} = \sum_{j=0}^p \frac{\mu_j}{j!}(it)^j + o(t^p)$$

and

$$l(t) := \log \mathbf{E}\{e^{it\xi}\} = \sum_{j=0}^{p} \frac{\kappa_j}{j!}(it)^j + o(t^p), \qquad t \to 0.$$

Here "log" denotes the principal branch of logarithm. The coefficients in the expansion of $f$ are the moments $\{\mu_j\}$ of $\xi$, and the coefficients in the expansion of $l$ are the cumulants $\{\kappa_j\}$ of $\xi$. For small $j$ one easily expresses $\kappa_j$ via $\mu_1, \mu_2, \dots, \mu_j$. In particular,

$$(3.4) \quad \kappa_1 = \mu_1, \qquad \kappa_2 = \mu_2 - \mu_1^2 = \mathbf{Var}\{\xi\}, \qquad \kappa_3 = \mu_3 - 3\mu_2\mu_1 + 2\mu_1^3,$$

$$\kappa_4 = \mu_4 - 3\mu_2^2 - 4\mu_3\mu_1 + 12\mu_2\mu_1^2 - 6\mu_1^4, \dots.$$

In general,

$$(3.5) \qquad \kappa_j = \sum_{\lambda} c_\lambda \mu_\lambda,$$

where the sum is over all additive partitions $\lambda$ of the set $\{1, \dots, j\}$, $c_\lambda$ are known coefficients and $\mu_\lambda = \prod_{l \in \lambda} \mu_l$; see, for example, [27].

We have [18]:

PROPOSITION 3.1. *Let $\xi$ be a random variable such that $\mathbf{E}\{|\xi|^{p+2}\} < \infty$ for a certain nonnegative integer $p$. Then for any function $\Phi : \mathbb{R} \to \mathbb{C}$ of the class $C^{p+1}$ with bounded derivatives $\Phi^{(l)}$, $l = 1, \dots, p+1$, we have*

$$(3.6) \qquad \mathbf{E}\{\xi\Phi(\xi)\} = \sum_{l=0}^{p} \frac{\kappa_{l+1}}{l!} \mathbf{E}\{\Phi^{(l)}(\xi)\} + \varepsilon_p,$$

*where the remainder term $\varepsilon_p$ admits the bound*

$$(3.7) \quad |\varepsilon_p| \le C_p \mathbf{E}\{|\xi|^{p+2}\} \sup_{t \in \mathbb{R}} |\Phi^{(p+1)}(t)|, \qquad C_p \le \frac{1 + (3+2p)^{p+2}}{(p+1)!}.$$

PROOF. Expanding the left- and the right-hand side of the identity

$$\mathbf{E}\{\xi e^{it\xi}\} = f(t)l'(t)$$

in powers of $it$, we obtain

$$(3.8) \qquad \mu_{r+1} = \sum_{j=0}^{r} \binom{r}{j} \kappa_{j+1}\mu_{r-j}, \qquad r = 0, 1, \dots, p.$$

Let $\pi$ be a polynomial of degree less or equal $p$. Then (3.8) implies that (3.6) is exact for $\Phi = \pi$, that is, is valid with $\varepsilon_p = 0$:

$$\mathbf{E}\{\xi\pi(\xi)\} = \sum_{j=0}^{p} \frac{\kappa_{j+1}}{j!}\mathbf{E}\{\pi^{(j)}(\xi)\}.$$

In the general case we write by Taylor's theorem $\Phi = \pi_p + r_p$, where $\pi_p$ is a polynomial of degree $p$ and

$$r_p(t) = \frac{t^{p+1}}{p!} \int_0^1 \Phi^{(p+1)}(tv)(1-v)^p \, dv.$$

Thus,

$$(3.9) \qquad |\mathbf{E}\{\xi\Phi(\xi)\} - \mathbf{E}\{\xi\pi_p(\xi)\}| \le \mathbf{E}\{|\xi r_p(\xi)|\} \le \frac{K_\Phi}{(p+1)!}\mathbf{E}\{|\xi|^{p+2}\},$$

where

$$K_\Phi = \sup_{t \in \mathbb{R}} |\Phi^{(p+1)}(t)| < \infty.$$

Besides,

$$\Phi^{(l)}(t) - \pi_p^{(l)}(t) = \frac{t^{p+1-l}}{(p-l)!} \int_0^1 \Phi^{(p+1)}(tv)(1-v)^{p-l} \, dv, \qquad l = 0, \dots, p,$$

and, therefore,

$$(3.10) \qquad \left| \mathbf{E}\{\xi\pi_p(\xi)\} - \sum_{j=0}^p \frac{\kappa_{j+1}}{j!}\mathbf{E}\{\Phi^{(j)}(\xi)\} \right| \le K_\Phi \sum_{j=0}^p \frac{|\kappa_{j+1}|\mathbf{E}\{|\xi|^{p-j+1}\}}{j!(p-j+1)!}.$$

The sum on the r.h.s. can be estimated with the help of the bound [27]:

$$(3.11) \qquad |\kappa_j| \le j^j \mathbf{E}\{|\xi - \mathbf{E}\{\xi\}|^j\}.$$

Since $(a+b)^j \le 2^{j-1}(a^j + b^j)$ for a positive integer $j$ and nonnegative $a$ and $b$, we have

$$(3.12) \qquad |\kappa_j| \le j^j \mathbf{E}\{(|\xi| + |\mathbf{E}\{\xi\}|)^j\} \le (2j)^j \mathbf{E}\{|\xi|^j\}.$$

This bound and the Hölder inequality $\mathbf{E}\{|\xi|^j\} \le \mathbf{E}\{|\xi|^{p+2}\}^{j/(p+2)}$ yield

$$(3.13) \qquad \begin{aligned} \sum_{j=0}^p \frac{|\kappa_{j+1}|\mathbf{E}\{|\xi|^{p-j+1}\}}{j!(p-j+1)!} &\le \mathbf{E}\{|\xi|^{p+2}\} \sum_{j=0}^p \frac{[2(j+1)]^{j+1}}{j!(p-j+1)!} \\ &\le \mathbf{E}\{|\xi|^{p+2}\}\frac{(3p+2)^{p+1}}{(p+1)!}. \end{aligned}$$

The proposition now follows from (3.9)–(3.13). $\quad\square$

Here is a simple "interpolation" corollary showing the mechanism of proximity of expectations with respect to the probability law of an arbitrary random variable and the Gaussian random variable with the same first and second moments. Its multivariate version will be often used below.

COROLLARY 3.1.   *Let $\xi$ be a random variable such that $\mathbf{E}_\xi\{|\xi|^{p+2}\} < \infty$ for a certain integer $p \geq 1$, $\mathbf{E}_\xi\{\xi\} = 0$, and let $\widehat{\xi}$ be the Gaussian random variable, whose first and second moments coincide with those of $\xi$. Then for any function $\Phi : \mathbb{R} \to \mathbb{C}$ of the class $C^{p+2}$ with bounded derivatives, we have*

$$\mathbf{E}_\xi\{\Phi(\xi)\} - \mathbf{E}_{\widehat{\xi}}\{\Phi(\widehat{\xi})\} = \sum_{l=2}^{p} \frac{\kappa_{l+1}}{2 l!} \int_0^1 \mathbf{E}\{\Phi^{(l+1)}(\xi(s))\} s^{(l-1)/2}\, ds + \varepsilon_p',$$

(3.14)

*where the symbols $\mathbf{E}_\xi\{\dots\}$ and $\mathbf{E}_{\widehat{\xi}}\{\dots\}$ denote the expectation with respect to the probability law of $\xi$ and $\widehat{\xi}$, $\{\kappa_j\}$ are the cumulants of $\xi$, $\mathbf{E}\{\dots\}$ denotes the expectation with respect to the product of probability laws of $\xi$ and $\widehat{\xi}$,*

(3.15)          $\xi(s) = s^{1/2}\xi + (1-s)^{1/2}\widehat{\xi}, \qquad 0 \leq s \leq 1,$

(3.16)          $|\varepsilon_p'| \leq C_p \mathbf{E}\{|\xi|^{p+2}\} \sup_{t \in \mathbb{R}} |\Phi^{(p+2)}(t)|$

*and $C_p$ satisfies (3.7).*

PROOF.   It suffices to write

$$\mathbf{E}_\xi\{\Phi(\xi)\} - \mathbf{E}_{\widehat{\xi}}\{\Phi(\widehat{\xi})\}$$

(3.17)      $$= \int_0^1 \frac{d}{ds} \mathbf{E}\{\Phi(\xi(s))\}\, ds$$

$$= \frac{1}{2} \int_0^1 \mathbf{E}\{s^{-1/2}\xi \Phi'(\xi(s)) - (1-s)^{-1/2}\widehat{\xi}\Phi'(\xi(s))\}\, ds$$

and use (3.6) for the first term in the parentheses and (2.20) for the second term.   $\square$

### 3.2. *Limiting normalized counting measure of eigenvalues.*   We will also need an analog of Theorem 2.1 on the limiting expectation of linear eigenvalue statistics of Wigner matrices. It has been known since the late 50s that the measure is again the semicircle law (2.31) (see [4, 12, 22] and [23] for results and references). We give below a new proof of this fact [25] that is based on the matrix analog of the "interpolation trick" (3.14) and illustrates the mechanism of coincidence of certain asymptotic results for Gaussian and non-Gaussian random matrices. The trick will be systematically used in what follows.

THEOREM 3.1.   *Let $M = n^{-1/2}W$ be the Wigner matrix (3.1)–(3.3), satisfying the condition*

(3.18)          $w_3 := \sup_n \max_{1 \leq j \leq k \leq n} \mathbf{E}\{|W_{jk}^{(n)}|^3\} < \infty$

*and $N_n$ be the normalized counting measure of its eigenvalues (2.33). Then*

$$\lim_{n \to \infty} \mathbf{E}\{N_n\} = N_{\mathrm{scl}},$$

*where $N_{\mathrm{scl}}$ is the semicircle law (2.31) and the convergence is understood as the weak convergence of measures.*

PROOF. It follows from (3.2) that we have (2.34) for the Wigner matrices. Thus, the sequence $\{\mathbf{E}\{N_n\}\}_{n \geq 0}$ is tight, and it suffices to prove its vague convergence, or, in view of the one-to-one correspondence between the nonnegative measures and their Stieltjes transforms (see, e.g., [1]), it suffices to prove the convergence of the Stieltjes transform of expectation of the normalized counting measure [see (2.35) and (2.37)] on a compact set of $\mathbb{C} \setminus \mathbb{R}$. Let $\widehat{M} = n^{-1/2}\widehat{W}$ be the GOE matrix (2.22)–(2.24), and $\widehat{G}(z)$ be its resolvent. Then, by Theorem 2.1, it suffices to prove that the limit

$$(3.19) \qquad \lim_{n \to \infty} |\mathbf{E}\{n^{-1}\operatorname{Tr} G(z)\} - \mathbf{E}\{n^{-1}\operatorname{Tr} \widehat{G}(z)\}| = 0$$

holds uniformly on a compact set of $\mathbb{C} \setminus \mathbb{R}$.

Following the idea of Corollary 3.1, consider the "interpolating" random matrix [cf. (3.15)]

$$(3.20) \qquad M(s) = s^{1/2}M + (1-s)^{1/2}\widehat{M}, \qquad 0 \leq s \leq 1,$$

viewed as defined on the product of the probability spaces of matrices $W$ and $\widehat{W}$. In other words, we assume that $W$ and $\widehat{W}$ in (3.20) are independent. We denote again by $\mathbf{E}\{\ldots\}$ the corresponding expectation in the product space. It is evident that $M(1) = M$, $M(0) = \widehat{M}$. Hence, if $G(s,z)$ is the resolvent of $M(s)$, then we have

$$\begin{aligned}
(3.21) \quad & n^{-1}\mathbf{E}\{\operatorname{Tr} G(z) - \operatorname{Tr} \widehat{G}(z)\} \\
& = \int_0^1 \frac{\partial}{\partial s}\mathbf{E}\{n^{-1}\operatorname{Tr} G(s,z)\}\,ds \\
& = -\frac{1}{2n^{3/2}}\int_0^1 \mathbf{E}\left\{\operatorname{Tr} \frac{\partial}{\partial z}G(s,z)(s^{-1/2}W - (1-s)^{-1/2}\widehat{W})\right\},
\end{aligned}$$

where we used (2.39) and (3.20).

Now we apply the differentiation formula (3.6) to transform the contribution of the first term in the parentheses of the r.h.s. of (3.21). To this end, we use the symmetry of the matrix $\{G_{jk}\}$ to write the corresponding expression as

$$(3.22) \qquad (n^3 s)^{-1/2}\sum_{1 \leq j \leq k \leq n} \beta_{jk}\mathbf{E}\{W_{jk}^{(n)}(G')_{jk}\},$$

where $\beta_{jk}$ are defined in (2.18) and we denote here and below

$$G' = \frac{\partial}{\partial z} G(s, z).$$

Since the random variables $W_{jk}^{(n)}$, $1 \leq j \leq k \leq n$, are independent, we can apply (3.6) with $p = 1$ and $\Phi = G'_{jk}$ to every term of the sum of (3.22). We obtain, in view of (3.2), (3.18) and (3.20),

$$(3.23) \quad \frac{w^2}{n^2} \sum_{1 \leq j \leq k \leq n} \mathbf{E}\{D_{jk}(s)(G')_{kj}\} + \varepsilon_1, \qquad D_{jk}(s) = \partial/\partial M_{jk}(s),$$

where [cf. (3.7)]

$$(3.24) \quad |\varepsilon_1| \leq \frac{C_1 w_3}{n^{5/2}} \sum_{1 \leq j \leq k \leq n} \sup_{M(s) \in \mathcal{S}_n} |D_{jk}^2(s)(G')_{jk}|,$$

$\mathcal{S}_n$ is the set of $n \times n$ real symmetric matrices, and $C_1$ is given by (3.7) for $p = 1$.

On the other hand, applying to the second term in the parentheses of (3.21) the Gaussian differential formula (2.20), we obtain again the first term of (3.23). Thus, the integrand of the r.h.s. of (3.21) is equal to $\varepsilon_1$.

It follows from (2.42) and its iterations that

$$(3.25) \quad |D_{jk}^l G_{jk}| \leq c_l/|\Im z|^{(l+1)}, \qquad |D_{jk}^l (G')_{jk}| \leq c_l/|\Im z|^{(l+2)},$$

where $c_l$ is an absolute constant for every $l$. The bounds and (3.24) imply

$$|\varepsilon_1| \leq \frac{C w_3}{n^{1/2} |\Im z|^4}$$

and $C$ denotes here and below a quantity that does not depend on $j$, $k$ and $n$, and may be distinct on different occasions.

This and the interpolation property (3.21) yield the assertion of the theorem. $\square$

In fact, we have more (see [4, 12] and [22] for other proofs and references).

THEOREM 3.2. *The assertion of Theorem 3.1 remains true under the condition*

$$(3.26) \quad \lim_{n \to \infty} n^{-2} \sum_{j,k=1}^n \int_{|W| > \tau \sqrt{n}} W^2 F_{jk}^{(n)}(dW) = 0 \qquad \forall \tau > 0.$$

PROOF.    Given $\tau > 0$, introduce the truncated matrix

$$M^\tau = \frac{W^\tau}{n^{1/2}}, \qquad W^\tau = \{W_{jk}^{(n)\tau} = \operatorname{sign} W_{jk}^{(n)} \max\{|W_{jk}^{(n)}|, \tau n^{1/2}\}\}_{j,k=1}^n.$$

(3.27)

Let $\mu_{l,jk}^\tau$ ($\mu_{l,jk}$) and $\kappa_{l,jk}^\tau$ ($\kappa_{l,jk}$) be the $l$th moment and cumulant of $W_{jk}^{(n)\tau}$ ($W_{jk}^{(n)}$), respectively. Then

$$|\mu_{l,jk}^\tau - \mu_{l,jk}| \le 2 \int_{|W|>\tau\sqrt{n}} |W|^l F_{jk}^{(n)}(dW).$$

This and (3.5) yield

$$(3.28) \qquad |\kappa_{l,jk}^\tau - \kappa_{l,jk}| \le C \int_{|W|>\tau\sqrt{n}} |W|^l F_{jk}^{(n)}(dW),$$

where $C$ depends only on $l$. In particular, we have

$$(3.29) \qquad |\kappa_{1,jk}^\tau - \kappa_{1,jk}| \le \frac{C}{\tau\sqrt{n}} \int_{|W|>\tau\sqrt{n}} W^2 F_{jk}^{(n)}(dW)$$

and

$$(3.30) \quad |\kappa_{l,jk}^\tau - \kappa_{l,jk}| \le \frac{C}{(\tau\sqrt{n})^{4-l}} \int_{|W|>\tau\sqrt{n}} W^4 F_{jk}^{(n)}(dW), \qquad l \le 4.$$

As in the previous theorem, it suffices to prove the limiting relation (3.19). To this end, we show first that, for every $\tau > 0$, the limit

$$(3.31) \qquad \lim_{n\to\infty} |\mathbf{E}\{n^{-1}\operatorname{Tr} G(z)\} - \mathbf{E}\{n^{-1}\operatorname{Tr} G^\tau(z)\}| = 0$$

with $G^\tau(z) = (M^\tau - z)^{-1}$ uniform on any compact set of $\mathbb{C} \setminus \mathbb{R}$. Indeed, we have by the resolvent identity (2.38) and the bound $|(G^\tau(z)G(z))_{jk}| \le |\Im z|^{-2}$,

$$\begin{aligned}
|\mathbf{E}\{n^{-1}&\operatorname{Tr}(G(z) - G^\tau(z))\}| \\
&= \left| \frac{1}{n^{3/2}} \sum_{j,k=1}^n \mathbf{E}\{(G^\tau(z)G(z))_{jk}(W_{jk}^{(n)} - W_{jk}^{(n)\tau})\} \right| \\
&\le \frac{1}{n^{3/2}|\Im z|^2} \sum_{j,k=1}^n \int_{|W|>\tau\sqrt{n}} |W| F_{jk}^{(n)}(dW) \le \frac{1}{|\Im z|^2\tau} L_n(\tau),
\end{aligned}$$

where

$$(3.32) \qquad L_n(\tau) = n^{-2} \sum_{j,k=1}^n \int_{|W|>\tau\sqrt{n}} W^2 F_{jk}^{(n)}(dW).$$

The last inequality and (3.26) imply (3.31).

Hence, it suffices to show that

$$(3.33) \qquad R_{n\tau} := \mathbf{E}\{n^{-1}\operatorname{Tr} G^{\tau}(z)\} - \mathbf{E}\{n^{-1}\operatorname{Tr} \widehat{G}(z)\}$$

vanishes after the subsequent limits

$$(3.34) \qquad n \to \infty, \qquad \tau \to 0,$$

uniformly on any compact set of $\mathbb{C} \setminus \mathbb{R}$.

Introduce the interpolation matrix

$$(3.35) \qquad M^{\tau}(s) = s^{1/2}M^{\tau} + (1-s)^{1/2}\widehat{M}, \qquad 0 \le s \le 1$$

[cf. (3.20)], denote its resolvent by $G(s, z) = (M^{\tau}(s) - z)^{-1}$, and get an analog of (3.21):

$$R_{n\tau} = -\frac{1}{2n^{3/2}} \int_0^1 \sum_{j,k=1}^n \mathbf{E}\{(G')_{jk}(s^{-1/2}W_{jk}^{(n)\tau} - (1-s)^{-1/2}\widehat{W}_{jk})\}\, ds.$$

As in the previous theorem, we apply the differentiation formula (3.6) with $p = 1$ to every term containing the factors $W_{jk}^{(n)\tau}$ and Gaussian differentiation formula (2.20) to every term containing $\widehat{W}_{jk}$, and obtain

$$\begin{aligned}
R_{n\tau} = -\frac{1}{2} \int_0^1 \Bigg( &\frac{1}{\sqrt{n^3 s}} \sum_{j,k=1}^n \kappa_{1,jk}^{\tau} \mathbf{E}\{(G')_{jk}\} \\
&+ \frac{1}{n^2} \sum_{j,k=1}^n \kappa_{2,jk}^{\tau} \mathbf{E}\{D_{jk}(s)(G')_{jk}\} \\
&+ \varepsilon_1 - \frac{1}{n^2} \sum_{j,k=1}^n w^2(1+\delta_{jk}) \mathbf{E}\{D_{jk}(s)(G')_{jk}\} \Bigg)\, ds,
\end{aligned}$$

where [cf. (3.24)]

$$|\varepsilon_1| \le \frac{C_1 s^{1/2}}{n^{5/2}} \sum_{1 \le j \le k \le n} \mathbf{E}\{|W_{jk}^{(n)\tau}|^3\} \sup_{M(s) \in \mathcal{S}_n} |D_{jk}^2(s)(G')_{jk}| \le \frac{2w^2 C_1 s^{1/2}}{|\Im z|^4}\tau,$$

and we took into account (3.25) and the bound

$$\mathbf{E}\{|W_{jk}^{(n)\tau}|^3\} \le \tau n^{1/2} w^2 (1+\delta_{jk}).$$

Besides, we have by (3.25), (3.29) and (3.28) with $l = 2$, and the equalities $\kappa_{1,jk} = 0$, $\kappa_{2,jk} = w^2(1+\delta_{jk})$,

$$\left| \frac{1}{n^{3/2}} \sum_{j,k=1}^n \kappa_{1,jk}^{\tau} \mathbf{E}\{(G')_{jk}\} \right| \le \frac{C}{\tau|\Im z|^2} L_n(\tau),$$

$$\left| \frac{1}{n^2} \sum_{j,k=1}^n (\kappa_{2,jk}^{\tau} - w^2(1+\delta_{jk})) \mathbf{E}\{D_{jk}(s)(G')_{jk}\} \right| \le \frac{C}{|\Im z|^3} L_n(\tau).$$

The last inequalities show that $R_{n\tau}$ of (3.33) vanishes after the subsequent limits (3.34). This completes the proof of the theorem. $\quad\square$

REMARK 3.1. Condition (3.26) is a matrix analog of the well-known Lindeberg condition

$$(3.36) \qquad \lim_{n\to\infty} \frac{1}{n} \sum_{j=1}^{n} \int_{|x|>\tau\sqrt{n}} x^2 F_j^{(n)}(dx) = 0 \qquad \forall \tau > 0$$

for a collection $\{\xi_j^{(n)}\}_{j=1}^{n}$ of independent random variables with probability laws $F_j^{(n)}$, $j = 1, \ldots, n$, to satisfy the central limit theorem. According to Theorem 3.2, the matrix analog (3.26) of the Lindeberg condition is sufficient for the validity of the semi-circle law for the Wigner ensembles. Thus, we can say that the semi-circle law is a universal limiting eigenvalues distribution of the Wigner ensembles in the same sense as the Gaussian distribution (the normal law) is universal for properly normalized sums of independent random variables. We mention two sufficient conditions for (3.26) to be valid, analogous to those of probability theory. The first is

$$\sup_n \max_{1 \le j \le k \le n} \mathbf{E}\{|W_{jk}^{(n)}|^{2+\delta}\} < \infty$$

for some $\delta > 0$. This is an analog of the Lyapunov condition of probability theory. The second sufficient condition requires that $\{W_{jk}^{(n)}\}_{1 \le j < k \le n}$ and $\{W_{jj}^{(n)}\}_{1 \le j \le n}$ are two collections of independent identically distributed random variables, whose probability laws $F_1$ and $F_2$ do not depend on $n$ and satisfy (3.2). This case generalizes the GOE, where $F_{1,2}$ are both Gaussian.

3.3. *Central limit theorem for linear eigenvalue statistics in the case of zero excess.* We consider here a particular case of the Wigner ensembles, for which the fourth cumulant of entries, known in statistics as the excess, is zero. The case is of interest because here the limiting Gaussian law has the same variance as in the GOE case (Theorem 2.2), moreover, it can be obtained from that for the GOE by applying the "interpolation" trick that was used in the proof of Theorems 3.1 and 3.2 [see formulas (3.20) and (3.35)]. We start from an analog of Theorem 3.1.

THEOREM 3.3. *Let $M = n^{-1/2}W$ be the real symmetric Wigner matrix (3.1)–(3.3). Assume the following:*

(i)

$$(3.37) \qquad w_5 := \sup_{n \in \mathbb{N}} \max_{1 \le j, k \le n} \mathbf{E}\{|W_{jk}^{(n)}|^5\} < \infty;$$

(ii) *the third and fourth moments do not depend on $j$, $k$ and $n$:*

$$(3.38) \qquad \mu_3 = \mathbf{E}\{(W_{jk}^{(n)})^3\}, \qquad \mu_4 = \mathbf{E}\{(W_{jk}^{(n)})^4\};$$

(iii) *the fourth cumulant of off-diagonal entries is zero:*

$$(3.39) \qquad \kappa_4 = \mu_4 - 3w^4 = 0.$$

*Let $\varphi \colon \mathbb{R} \to \mathbb{R}$ be a test function whose Fourier transform $\widehat{\varphi}$ (2.53) satisfies the condition*

$$(3.40) \qquad \int (1 + |t|^5)|\widehat{\varphi}(t)|\, dt < \infty.$$

*Then the corresponding centered linear eigenvalue statistic $\mathcal{N}_n^\circ[\varphi]$ [see (2.47)] converges in distribution to the Gaussian random variable of zero mean and variance $V_{\mathrm{GOE}}[\varphi]$ of (2.48).*

REMARK 3.2.   It may seem not too natural to have the $(j, k)$ dependent second moments (3.2) of $W_{jk}^{(n)}$ and the $(j, k)$ independent fourth moments (3.38). This is only for the sake of technical simplicity of the proof. In fact, it can be shown that the result does not depend on the diagonal entries, in particular, we can assume that the second moments will be the same for all $1 \le j \le k \le n$, or that only the fourth moments of off-diagonal entries are the same and the fourth moments of diagonal entries are just uniformly bounded. Likewise, we can replace (3.39) by

$$\lim_{n \to \infty} \sup_n \max_{1 \le j < k \le n} |\kappa_{4,jk}| = 0,$$

where $\kappa_{4,jk}$ is the fourth cumulants of $W_{jk}^{(n)}$.

PROOF OF THEOREM 3.3.   Let $\widehat{M} = n^{-1/2}\widehat{W}$ be the GOE matrix (2.22)–(2.24) with the same variance of entries as the Wigner matrix, and $\widehat{\mathcal{N}}_n^\circ[\varphi]$ be the centered linear eigenvalue statistic of the GOE. Then, in view of Theorem 2.2, it suffices to show that, for every $x \in \mathbb{R}$,

$$(3.41) \qquad R_n(x) := \mathbf{E}\{e^{ix\mathcal{N}_n^\circ[\varphi]}\} - \mathbf{E}\{e^{ix\widehat{\mathcal{N}}_n^\circ[\varphi]}\} \to 0, \qquad n \to \infty.$$

Denoting

$$(3.42) \qquad e_n(s, x) = \exp\{ix\,\mathrm{Tr}\,\varphi(M(s))^\circ\},$$

where $M(s)$ is the interpolating matrix (3.20), we have

$$R_n(x) = \int_0^1 \frac{\partial}{\partial s}\mathbf{E}\{e_n(s, x)\}\, ds$$

$$= \frac{ix}{2\sqrt{n}} \int_0^1 \mathbf{E}\{e_n^\circ(s, x)\,\mathrm{Tr}\,\varphi'(M(s))(s^{-1/2}W - (1-s)^{-1/2}\widehat{W})\}\, ds$$

[cf. (3.17)] or, after using (2.73),

$$R_n(x) = -\frac{x}{2} \int_0^1 ds \int t\widehat{\varphi}(t)[A_n - B_n]\, dt, \tag{3.43}$$

where

$$A_n = \frac{1}{\sqrt{ns}} \sum_{j,k=1}^n \mathbf{E}\{W_{jk}^{(n)}\Phi_n\}, \tag{3.44}$$

$$B_n = \frac{1}{\sqrt{n(1-s)}} \sum_{j,k=1}^n \mathbf{E}\{\widehat{W}_{jk}\Phi_n\} \tag{3.45}$$

with

$$\Phi_n = e_n^\circ(s,x)U_{jk}(s,t), \qquad U(s,t) = e^{itM(s)}. \tag{3.46}$$

Applying (3.6) with $p=3$ to every term of (3.44) and (2.20) to every term of (3.45), we obtain (cf. Corollary 3.1)

$$A_n - B_n = T_2 + T_3 + \varepsilon_3, \tag{3.47}$$

where [cf. (3.23)]

$$T_l = \frac{s^{(l-1)/2}}{l!n^{(l+1)/2}} \sum_{j,k=1}^n \kappa_{l+1,jk}\mathbf{E}\{D_{jk}^l(s)\Phi_n\}, \qquad D_{jk}(s) = \partial/\partial M_{jk}(s), \tag{3.48}$$

and by (3.7) and (3.37) [cf. (3.24)],

$$|\varepsilon_3| \leq \frac{C_3 w_5}{n^{5/2}} \sum_{j,k=1}^n \sup_{M\in\mathcal{S}_n} |D_{jk}^4(s)\Phi_n|_{M(s)=M}|. \tag{3.49}$$

It follows then from (2.17), (2.19), (2.72) and (3.40) that

$$|D_{jk}^l(s)\Phi_n| \leq C_l(t,x), \qquad 0 \leq l \leq 4, \tag{3.50}$$

and we denote here and below $C_l(t,x)$ a polynomial in $|t|$ and $|x|$ of degree $l$, independent of $j,k$ and $n$ and not necessarily the same at each occurrence. This implies that

$$|\varepsilon_3| \leq C_4(t,x)/n^{1/2}. \tag{3.51}$$

Furthermore, it follows from (3.2) and (3.39) that $\kappa_{4,jk} = -9\delta_{jk}w^4$. This, (3.48) for $l=3$, and (3.50) yield

$$|T_3| \leq C_3(t,x)/n. \tag{3.52}$$

To get a vanishing bound for $T_2$ of (3.48), we use (2.17) and (2.72) to find the second derivative $D^2_{jk}\Phi_n$ and we take into account (3.38) to have

$$
\begin{aligned}
T_2 = -\frac{\sqrt{s}\mu_3}{n^{3/2}} \sum_{j,k=1}^n \beta_{jk}^2 \mathbf{E}\Big\{ & e_n^\circ [(U_{jk} * U_{jk} * U_{jk})(t) + 3(U_{jk} * U_{jj} * U_{kk})(t)] \\
& + 2xe_n[(U_{jk} * U_{jk})(t) + (U_{jj} * U_{kk})(t)] \\
& \times \int \theta \widehat{\varphi}(\theta) U_{jk}(\theta)\, d\theta \\
& - 2x^2 e_n U_{jk}(t) \Big(\int \theta \widehat{\varphi}(\theta) U_{jk}(\theta)\, d\theta\Big)^2 \\
& + ixe_n U_{jk}(t) \int \theta \widehat{\varphi}(\theta)[(U_{jk} * U_{jk})(\theta) \\
& \hspace{6cm} - (U_{jj} * U_{kk})(\theta)]\, d\theta \Big\},
\end{aligned}
$$
(3.53)

where we write $e_n$ for $e_n(s,x)$, $U(t)$ for $U(s,t)$ and take into account that the convolution operation "$*$" of Proposition 2.1(iii) is commutative.

Consider the two types of the sums above:

$$
\begin{aligned}
T_{21} &= n^{-3/2} \sum_{j,k=1}^n U_{jk}(t_1) U_{jj}(t_2) U_{kk}(t_3), \\
T_{22} &= n^{-3/2} \sum_{j,k=1}^n U_{jk}(t_1) U_{jk}(t_2) U_{jk}(t_3).
\end{aligned}
$$
(3.54)

It follows from the Schwarz inequality and (2.16) that

$$
\sum_{j,k=1}^n |U_{jk}(t_1) U_{jk}(t_2)| \le \Big(\sum_{j,k=1}^n |U_{jk}(t_1)|^2\Big)^{1/2} \Big(\sum_{j,k=1}^n |U_{jk}(t_2)|^2\Big)^{1/2} = n,
$$

hence, $|T_{22}| \le n^{-1/2}$. Besides, writing

$$
T_{21} = n^{-1/2}(U(t_1)V(t_2), V(t_3)), \qquad V(t) = n^{-1/2}(U_{11}(t), \dots, U_{nn}(t))^T,
$$

where by (2.16) $\|V(t)\| \le 1$, $\|U(t)\| = 1$, we conclude that $|T_{21}| \le n^{-1/2}$, hence,

$$
|T_2| \le C_2(t,x)/n^{1/2}.
$$
(3.55)

This together with (3.40), (3.47), (3.51) and (3.52) imply that the r.h.s. of (3.43) is $O(n^{-1/2})$ as $n \to \infty$. We obtain (3.41), hence the assertion of the theorem. $\square$

In fact, we have more (see [18] for a particular case of traces of resolvent).

THEOREM 3.4.    *Theorem 3.3 remains valid if its condition (3.37) is re-
placed by the Lindeberg type condition for the fourth moments of entries of
W*

$$\lim_{n \to \infty} L_n^{(4)}(\tau) = 0 \qquad \forall \tau > 0, \tag{3.56}$$

*where [cf. (3.32)]*

$$L_n^{(4)}(\tau) = \frac{1}{n^2} \sum_{j,k=1}^n \int_{|W| > \tau \sqrt{n}} W^4 F_{jk}^{(n)}(dW). \tag{3.57}$$

PROOF.    Consider again the truncated matrix $M^\tau$ of (3.27). Since

$$\mathbf{P}\{W \neq W^\tau\} \leq \sum_{j,k=1}^n \mathbf{P}\{W_{jk}^{(n)} \neq W_{jk}^{(n)\tau}\}$$

$$= \sum_{j,k=1}^n \int_{|W| > \tau \sqrt{n}} F_{jk}^{(n)}(dW) \leq \tau^{-4} L_n^{(4)}(\tau), \tag{3.58}$$

we have, in view of (3.56),

$$\lim_{n \to \infty} \mathbf{E}\{e^{ix\mathcal{N}_n^\circ[\varphi]} - e^{ix\mathcal{N}_{n\tau}^\circ[\varphi]}\} = 0, \tag{3.59}$$

where $\mathcal{N}_{n\tau}[\varphi] = \mathrm{Tr}\,\varphi(M^\tau)$. Now it suffices to prove that if $\widehat{\mathcal{N}}_n[\varphi] = \mathrm{Tr}\,\varphi(\widehat{M})$
is the linear eigenvalue statistics of the GOE matrix $\widehat{M}$, then [cf. (3.41)]

$$R_n^\tau(x) = \mathbf{E}\{e^{ix\mathcal{N}_{n\tau}^\circ[\varphi]}\} - \mathbf{E}\{e^{ix\widehat{\mathcal{N}}_n^\circ[\varphi]}\} \tag{3.60}$$

vanishes after the limit (3.34). To this end, we use again the interpolation
matrix $M^\tau(s)$ of (3.35), and get analogs of (3.43)–(3.46):

$$R_n^\tau(x) = -\frac{x}{2} \int_0^1 ds \int t\widehat{\varphi}(t)[A_{n\tau} - B_{n\tau}]\,dt,$$

$$A_{n\tau} = \frac{1}{\sqrt{ns}} \sum_{j,k=1}^n \mathbf{E}\{W_{jk}^{(n)\tau} \Phi_{n\tau}\}, \tag{3.61}$$

$$B_{n\tau} = \frac{1}{\sqrt{n(1-s)}} \sum_{j,k=1}^n \mathbf{E}\{\widehat{W}_{jk} \Phi_{n\tau}\},$$

where now

$$\Phi_{n\tau} = e_{n\tau}^\circ(s,x)U_{jk}^\tau(s,t),$$

$$U^\tau(s,t) = \exp\{itM^\tau(s)\}, \qquad e_{n\tau}(s,x) = \exp\{ix\,\mathrm{Tr}\,\varphi(M^\tau(s))^\circ\}. \tag{3.62}$$

Applying (3.6) with $p = 3$, we obtain [cf. (3.47)–(3.49)]

$$(3.63) \quad A_{n\tau} = \sum_{l=0}^{3} T_{l\tau} + \varepsilon_{3\tau},$$

$$(3.64) \quad T_{l\tau} = \frac{s^{(l-1)/2}}{l! n^{(l+1)/2}} \sum_{j,k=1}^{n} \kappa_{l+1,jk}^{\tau} \mathbf{E}\{D_{jk}^{l}(s)\Phi_{n\tau}\}, \qquad l = 0, 1, 2, 3,$$

and

$$|\varepsilon_{3\tau}| \leq \frac{C_3}{n^{5/2}} \sum_{j,k=1}^{n} \mathbf{E}\{|W_{jk}^{(n)\tau}|^5\} \sup_{M \in \mathcal{S}_n} |D_{jk}^4(s)\Phi_{n\tau}|_{M(s)=M}|.$$

Since $\mathbf{E}\{|W_{jk}^{(n)\tau}|^5\} \leq \tau \sqrt{n}\mu_4$, then we have, in view of (3.50),

$$(3.65) \qquad\qquad |\varepsilon_{3\tau}| \leq C_4(t,x)\tau.$$

Besides, it follows from (3.50) and (3.30) that we can replace $T_{l\tau}$ of (3.64) by $T_l$ of (3.48) with $\Phi_{n\tau}$ of (3.62):

$$(3.66) \qquad\qquad T_{l\tau} = T_l + r_l,$$

where the error term $r_l$ satisfies

$$
(3.67)
\begin{aligned}
|r_l| &\leq \frac{s^{(l-1)/2}}{l! n^{(l+1)/2}} \sum_{j,k=1}^{n} |\kappa_{l+1,jk}^{\tau} - \kappa_{l+1,jk}| |\mathbf{E}\{D_{jk}^{l}(s)\Phi_{n\tau}\}| \\
&\leq s^{(l-1)/2} C_l(t,x) \tau^{l-3} L_n^{(4)}(\tau).
\end{aligned}
$$

We have by (3.2) $T_0 = 0$, $T_1 = B_n$, and $T_2$ and $T_3$ satisfy (3.55) and (3.52), respectively. This together with (3.40), (3.43) and (3.65) imply (3.60) and complete the proof of the theorem.  $\square$

3.4. *Central limit theorem in general case.* Here we prove the CLT for linear eigenvalue statistics of the Wigner random matrix not assuming that the fourth cumulant of its entries is zero [see (3.39)]. We use the scheme of the proof of Theorem 2.2, based on the Gaussian differentiation formula (2.20) and the Poincaré type "a priory" bound (2.26) for the variance of statistics. We have the extension of (2.20), given by (3.6). As for an analog of (2.26), it is given by the theorem below [see also (3.90)].

THEOREM 3.5. *Let $M = n^{-1/2}W$ be the Wigner matrix (3.1)–(3.3) satisfying (3.38) and (3.56), $M^{\tau}$ be corresponding truncated matrix (3.27), and*

$$(3.68) \qquad u_{n\tau}(t) = \operatorname{Tr} U^{\tau}(t), \qquad U^{\tau}(t) = \exp(it M^{\tau}).$$

*Then for any fixed $\tau > 0$,*

$$(3.69) \qquad \mathbf{Var}\{u_{n\tau}(t)\} \leq C_\tau(\mu_4)(1 + |t|^4)^2,$$

$$(3.70) \qquad \mathbf{Var}\{\mathcal{N}_{n\tau}[\varphi]\} \leq C_\tau(\mu_4)\left(\int (1 + |t|^4)|\widehat{\varphi}(t)|\, dt\right)^2,$$

*where $C_\tau(\mu_4)$ depends only on $\mu_4$ and $\tau$.*

PROOF.  Note first that by the Schwarz inequality for expectations and (2.60) we have

$$(3.71) \begin{aligned} \mathbf{Var}\{\mathcal{N}_{n\tau}[\varphi]\} &= \int\int \widehat{\varphi}(t_1)\overline{\widehat{\varphi}(t_2)}\mathbf{E}\{u_{n\tau}^\circ(t_1)\overline{u_{n\tau}^\circ(t_2)}\}\, dt_1\, dt_2 \\ &\leq \left(\int \mathbf{Var}^{1/2}\{u_{n\tau}(t)\}|\widehat{\varphi}(t)|\, dt\right)^2 \end{aligned}$$

and it suffices to get bound (3.69) for

$$V_n = \mathbf{Var}\{u_{n\tau}(t)\}.$$

Denoting $\widehat{u}_n(t) = \mathrm{Tr}\exp\{it\widehat{M}\}$, where $\widehat{M}$ is the GOE matrix, we write

$$(3.72) \begin{aligned} V_n &= \mathbf{E}\{\widehat{u}_n(t)u_{n\tau}^\circ(-t)\} + \mathbf{E}\{(u_{n\tau}(t) - \widehat{u}_n(t))\widehat{u}_n^\circ(-t)\} \\ &\quad + \mathbf{E}\{(u_{n\tau}(t) - \widehat{u}_n(t))(u_{n\tau}^\circ(-t) - \widehat{u}_n^\circ(-t))\} = K_1 + K_2 + K_3. \end{aligned}$$

We have, by the Schwarz inequality and (2.65),

$$(3.73) \qquad |K_1| \leq \sqrt{2}w|t|V_n^{1/2}, \qquad |K_2| \leq \sqrt{2}w|t|V_n^{1/2} + 2w^2t^2.$$

To estimate $K_3$, we use the interpolating matrix (3.35) to write

$$(3.74) \qquad K_3 = \frac{it}{2}\int_0^1 [A_n' - B_n']\, ds,$$

where

$$(3.75) \quad A_n' = \frac{1}{\sqrt{ns}}\sum_{j,k=1}^n \mathbf{E}\{W_{jk}^{(n)\tau}\Phi_n'\}, \qquad B_n' = \frac{1}{\sqrt{n(1-s)}}\sum_{j,k=1}^n \mathbf{E}\{\widehat{W}_{jk}\Phi_n'\}$$

with

$$(3.76) \qquad \Phi_n' = U_{jk}^\tau(s,t)(u_{n\tau}^\circ(-t) - \widehat{u}_n^\circ(-t))$$

and $U^\tau(s,t)$ being defined in (3.62).

Applying (3.6) with $p = 2$ and

$$\Phi(W) := \mathbf{E}\{\Phi_n'|_{M_{jk}^\tau(s) = (s/n)^{1/2}W + (1-s)^{1/2}\widehat{M}_{jk}}\}$$

to every term

$$\mathbf{E}\{W_{jk}^{(n)\tau}\Phi_n'\} = \int \Phi(W)W F_{jk}^{(n)\tau}(dW)$$

of the sum in $A_n'$ (3.75), we obtain

$$(3.77) \qquad A_n' = \sum_{l=0}^{2} T_{l\tau}' + \varepsilon_{2\tau},$$

where $T_{l\tau}'$ is defined by (3.64) with $\Phi_n'$ of (3.76) instead of $\Phi_n(s)$ of (3.62), and

$$(3.78) \qquad |\varepsilon_{2\tau}| \le \frac{C_2\mu_4}{n^2} \sum_{j,k=1}^{n} \sup_{|W| \le \tau\sqrt{n}} |\mathbf{E}\{D_{jk}^3(s)\Phi(W)\}|.$$

Since

$$(3.79) \qquad \begin{aligned} &\mathbf{E}\{D_{jk}^l(s)\Phi_n'\} \\ &= \mathbf{E}\{(u_{n\tau}^\circ(-t) - \widehat{u}_n^\circ(-t))D_{jk}^l(s)U_{jk}^\tau(s,t)\} \\ &\quad + \sum_{q=1}^{l}\binom{l}{q}\mathbf{E}\{D_{jk}^q(s)(u_{n\tau}(-t) - \widehat{u}_n(-t))D_{jk}^{l-q}(s)U_{jk}^\tau(s,t)\} \end{aligned}$$

and by (2.28),

$$(3.80) \qquad \frac{\partial}{\partial M_{jk}} u_{n\tau}(t) = 2i\beta_{jk}t U_{jk}^\tau(t),$$

the Schwarz inequality and (2.19) yield

$$(3.81) \qquad |\mathbf{E}\{D_{jk}^l(s)\Phi_n'\}| \le C_l(t)(V_n^{1/2} + 1).$$

Here and below we denote by $C_l(t)$ an $n$-independent polynomial in $|t|$ of degree $l$. This and (3.30) imply [cf. (3.66) and (3.67)]

$$(3.82) \qquad T_{l\tau}' = T_l' + r_l', \qquad l = 0,1,2,$$

where $T_l'$ is defined by (3.48) with $\Phi_n'$ of (3.76) instead of $\Phi_n$ of (3.46), $\kappa_{1,jk} = 0$, $\kappa_{2,jk} = (1 + \delta_{jk})w^2$, $\kappa_{3,jk} = \mu_3$, and

$$|r_l'| \le s^{(l-1)/2}C_l(t)\tau^{l-3}L_n^{(4)}(\tau)(V_n^{1/2} + 1).$$

Taking in account (3.56), we have for sufficiently large $n$

$$(3.83) \qquad |r_l'| \le s^{(l-1)/2}C_l(t)\tau^{l-3}(V_n^{1/2} + 1).$$

We see that $T_0' = 0$, and by applying (2.20) to $B_n'$ of (3.75), we have $T_1' = B_n'$. Besides, since by (3.79)

$$T_2' = \frac{s^{1/2}\mu_3}{2}\left(\mathbf{E}\left\{(u_{n\tau}^\circ(-t) - \widehat{u}_n^\circ(-t))\left(n^{-3/2}\sum_{j,k=1}^n D_{jk}^2(s)U_{jk}^\tau(s,t)\right)\right\}\right.$$

$$+ \sum_{q=1}^2 \binom{2}{q}\mathbf{E}\left\{n^{-3/2}\sum_{j,k=1} D_{jk}^q(s)(u_{n\tau}(-t) - \widehat{u}_n(-t))\right.$$

$$\left.\left.\times D_{jk}^{2-q}(s)U_{jk}^\tau(s,t)\right\}\right),$$

then using the Schwarz inequality and (2.19) to estimate the first term, and (3.80) and the argument leading to (3.55) to estimate the second term, we obtain

$$|T_2'| \le C_2(t)(V_n^{1/2} + 1).$$

It follows from the above for the integrand in (3.74)

$$(3.84) \qquad |A_n' - B_n'| \le |\varepsilon_{2\tau}| + C_2(t)(V_n^{1/2} + 1),$$

where $\varepsilon_{2\tau}$ is defined in (3.78) and we have, in view of (3.79),

$$(3.85) \qquad |\varepsilon_{2\tau}| \le \frac{C_2\mu_4}{n^2}\sum_{j,k=1}^n S_{jk} + |\varepsilon_{2\tau}'|,$$

where

$$S_{jk} = \sup_{|W| \le \tau\sqrt{n}}|\mathbf{E}\{(u_{n\tau}^\circ(-t) - \widehat{u}_n^\circ(-t))$$

$$\times D_{jk}^3 U_{jk}^\tau(s,t)|_{M_{jk}^\tau(s)=(s/n)^{1/2}W + (1-s)^{1/2}\widehat{M}_{jk}}\}|$$

and by (2.16), (2.19) and (3.80),

$$(3.86) \qquad |\varepsilon_{2\tau}'| \le C_3(t)$$

with $C_3(t)$ of (3.81).

To estimate $S_{jk}$, we repeat again the above interpolating procedure, and obtain for every fixed pair $\{j,k\}$

$$S_{jk} = \frac{|t|}{2}\sup_{|W| \le \tau\sqrt{n}}\left|\int_0^1 ds_1 \frac{1}{\sqrt{n}}\sum_{p,q=1}^n \mathbf{E}\{(s_1^{-1/2}W_{pq}^{(n)\tau} - (1-s_1)^{-1/2}\widehat{W}_{pq})\right.$$

$$\left.\times \Phi_n''|_{M_{jk}^\tau(s)=(s/n)^{1/2}W + (1-s)^{1/2}\widehat{M}_{jk}}\}\right|,$$

where

$$(3.87) \qquad \Phi_n'' = U_{pq}^\tau(s_1, t) D_{jk}^3 U_{jk}^\tau(s, t), \qquad |\Phi_n''| \le C_3(t).$$

By the condition $|W| \le \tau \sqrt{n}$ and (3.87), two terms of the sum corresponding to $W_{pq} = W_{jk} = W$ are bounded by $C_3(t)$ for every fixed $\tau > 0$. Hence, applying (3.6) and (2.20) to the rest of the terms, and using the notation $\sum_{p,q}'$ for the sum with $\{p, q\} \ne \{j, k\}$ and $\{p, q\} \ne \{k, j\}$, we obtain

$$(3.88) \qquad S_{jk} \le C_4(t) + \frac{|t|}{2} \sup_{|W| \le \tau \sqrt{n}} \int_0^1 |A_n'' - B_n''| \, ds_1$$

with

$$A_n'' = \frac{1}{\sqrt{s_1 n}} \sum_{p,q}' \mathbf{E}\{W_{pq}^{(n)\tau} \Phi_n''|_{M_{jk}^\tau(s)=(s/n)^{1/2}W+(1-s)^{1/2}\widehat{M}_{jk}}\} = \sum_0^2 T_{l\tau}'' + \varepsilon_{2\tau}''$$

and

$$B_n'' = \frac{w^2}{n} \sum_{p,q}' (1 + \delta_{pq}) \mathbf{E}\{D_{pq}(s_1) \Phi_n''|_{M_{jk}^\tau(s)=(s/n)^{1/2}W+(1-s)^{1/2}\widehat{M}_{jk}}\},$$

where [cf. (3.64)]

$$T_{l\tau}'' = \frac{s_1^{(l-1)/2}}{l! n^{(l+1)/2}} \sum_{p,q}' \kappa_{l+1,pq}^\tau \mathbf{E}\{D_{pq}^l(s_1) \Phi_n''|_{M_{jk}^\tau(s)=(s/n)^{1/2}W+(1-s)^{1/2}\widehat{M}_{jk}}\}$$

and

$$|\varepsilon_{2\tau}''| \le \frac{C_2 \mu_4}{n^2} \sum_{p,q}' \sup_{M \in S_n} |D_{pq}^3(s_1) \Phi_n''|_{M_{jk}^\tau(s)=(s/n)^{1/2}W+(1-s)^{1/2}\widehat{M}_{jk}, M^\tau(s_1)=M}|.$$

Since $|D_{pq}^l(s)\Phi_n''| \le C_{l+3}(t)$, then $|\varepsilon_{2\tau}''| \le C_6(t)$. Besides, in view of (3.30), we have an analog of (3.82) and (3.83):

$$T_{l\tau}'' = T_l'' + r_l'', \qquad l = 0, 1, 2,$$

where an argument, similar to that leading to (3.55), implies

$$|T_2''| \le C_5(t) n^{-1/2}.$$

We conclude that, for every $\tau > 0$,

$$\sup_{|W| \le \tau \sqrt{n}} \int_0^1 |A_n'' - B_n''| \, ds_1 \le C_6(t).$$

Plugging this estimate in (3.88) and then in (3.85), we obtain in view of (3.86) that

$$|\varepsilon_{2\tau}| \le C_7(t).$$

This, (3.74) and (3.84) imply inequality $|K_3| \leq C_3(t)V_n^{1/2} + C_8(t)$, which together with (3.72) and (3.73) allow us to write the quadratic inequality for $V_n$:

$$V_n \leq C_3(t)V_n^{1/2} + C_8(t)$$

valid for every fixed $\tau > 0$ and any real $t$ and implying (3.69). $\quad\square$

REMARK 3.3. A similar but much simpler argument allows us to prove that if

$$(3.89) \qquad w_6 := \sup_n \max_{1 \leq j < k \leq n} E\{|W_{jk}^{(n)}|^6\} < \infty,$$

then we have the bounds

$$(3.90) \qquad \mathbf{Var}\{\mathcal{N}_n[\varphi]\} \leq C(w_6)\left(\int (1 + |t|^3)|\widehat{\varphi}(t)|\, dt\right)^2$$

and

$$(3.91) \qquad \mathbf{Var}\{u_n(t)\} \leq C(w_6)(1 + |t|^3)^2,$$

where $C(w_6)$ depends only on $w_6$. The proof is based on the representation

$$\mathbf{Var}\{u_n(t)\} = \mathbf{E}\{\widehat{u}_n(t)u_n^\circ(-t)\} + \mathbf{E}\{(u_n(t) - \widehat{u}_n(t))u_n^\circ(-t)\},$$

the interpolation procedure, and the differentiation formula (3.6) with $p = 4$ in the second term.

Now we can prove the corresponding CLT.

THEOREM 3.6. Let $M = n^{-1/2}W$ be the real symmetric Wigner matrix (3.1)–(3.3), satisfying (3.38) and (3.56), and $\varphi : \mathbb{R} \to \mathbb{R}$ be a test function whose Fourier transform $\widehat{\varphi}$ satisfies (3.40). Then the centered linear eigenvalue statistic $\mathcal{N}_n^\circ[\varphi]$ [see (2.47)] converges in distribution to the Gaussian random variable of zero mean and variance

$$(3.92) \qquad V_{\mathrm{Wig}}[\varphi] = V_{\mathrm{GOE}}[\varphi] + \frac{\kappa_4}{2\pi^2 w^8}\left(\int_{-2w}^{2w} \varphi(\mu)\frac{2w^2 - \mu^2}{\sqrt{4w^2 - \mu^2}}\, d\mu\right)^2,$$

where $V_{\mathrm{GOE}}[\varphi]$ is given by (2.48), and $\kappa_4 = \mu_4 - 3w^4$ is the fourth cumulant of the off-diagonal entries of $W$.

PROOF. Following the scheme of the proof of Theorems 2.2, we show that the limit $Z(x)$ of characteristic functions $Z_n(x) = \mathbf{E}\{\exp(ix\mathcal{N}_n^\circ[\varphi])\}$ satisfies (2.55) with $V_{\mathrm{GOE}}$ of (2.48) replaced by $V_{\mathrm{Wig}}$ of (3.92). In view of (3.59), it suffices to find the limit as $n \to \infty$ of the characteristic functions

$$(3.93) \qquad Z_{n\tau}(x) = \mathbf{E}\{e_{n\tau}(x)\}, \qquad e_{n\tau}(x) = \exp\{ix\mathcal{N}_{n\tau}^\circ[\varphi]\}$$

of the centered statistics $\mathcal{N}_{n\tau}^{\circ}[\varphi]$ of truncated matrix $M^{\tau} = n^{-1/2}W^{\tau}$ of (3.27), and then pass to the limit $\tau \to 0$.

It is easy to see that formulas (2.57)–(2.64) with $u_{n\tau}$ and $Y_{n\tau}(x,t) = \mathbf{E}\{u_{n\tau}(t) \times e_{n\tau}^{\circ}(x)\}$ instead of $u_n$ and $Y_n(x,t)$ of (2.61) and (2.63) are valid in the Wigner case as well, and that (3.69) and (3.70) imply the analogs of (2.67) and (2.69) for $Y_{n\tau}$:

$$(3.94) \qquad |Y_{n\tau}(x,t)| \leq C_{\tau}^{1/2}(\mu_4)(1+|t|)^4$$

and

$$(3.95) \qquad \left|\frac{\partial}{\partial x}Y_{n\tau}(x,t)\right| \leq C_{\tau}(\mu_4)\left(\int(1+|t|^4)|\widehat{\varphi}(t)|\,dt\right)^2,$$

where $C_{\tau}(\mu_4)$ depends only on $\tau$ and $\mu_4$.

To prove the uniform boundedness of $\partial Y_n(x,t)/\partial t$ [an analog of (2.68)], we note first that, by (3.1) and (3.68),

$$(3.96) \qquad \frac{\partial}{\partial t}Y_{n\tau}(x,t) = \mathbf{E}\{u_{n\tau}'(t)e_{n\tau}^{\circ}(x)\} = \frac{i}{\sqrt{n}}\sum_{j,k=1}^{n}\mathbf{E}\{W_{jk}^{(n)\tau}\Phi_n\},$$

where

$$(3.97) \qquad \Phi_n = U_{jk}^{\tau}(t)e_{n\tau}^{\circ}(x), \qquad |D_{jk}^{l}\Phi_n| \leq C_l(t,x), \qquad 0 \leq l \leq 5$$

[see (3.50)]. Treating the r.h.s. of (3.96) as $A_{n\tau}$ of (3.63) and applying (3.6) with $p = 2$, we obtain

$$(3.98) \qquad \frac{\partial}{\partial t}Y_{n\tau}(x,t) = \frac{iw^2}{n}\sum_{j,k=1}^{n}(1+\delta_{jk})\mathbf{E}\{D_{jk}\Phi_n\} + O(1),$$

where the error term is bounded by $C_3(t,x)$ as $n \to \infty$ in view of (3.7), (3.55), (3.67) and (3.97). By using (2.17) and (2.72), we obtain for the first term of the r.h.s. of (3.98)

$$itw^2 n^{-1}Y_{n\tau}(x,t) + iw^2\int_0^t \mathbf{E}\{n^{-1}u_{n\tau}(t-t_1)\}Y_{n\tau}(x,t_1)\,dt_1$$

$$+ iw^2\int_0^t \mathbf{E}\{n^{-1}u_{n\tau}(t-t_1)u_{n\tau}^{\circ}(t_1)e_{n\tau}^{\circ}(x)\}\,dt_1$$

$$- 2w^2 x\int t_1\widehat{\varphi}(t_1)\mathbf{E}\{n^{-1}u_{n\tau}(t+t_1)e_{n\tau}(x)\}\,dt_1,$$

where the first two terms are bounded in view of (3.94), the last term is bounded by $2w^2|x|\int|t_1||\widehat{\varphi}(t_1)|\,dt_1$, and the third term satisfies

$$(3.99) \qquad |\mathbf{E}\{n^{-1}u_{n\tau}(t-t_1)u_{n\tau}^{\circ}(t_1)e_{n\tau}^{\circ}(x)\}| \leq 2C_{\tau}(\mu_4)^{1/2}(1+|t_1|)^4$$

in view of (3.69). It follows then from (3.96)–(3.99) that, for any fixed $\tau > 0$,

$$(3.100) \qquad \left| \frac{\partial}{\partial t} Y_{n\tau}(t,x) \right| \leq C_5(t,x).$$

Thus, we have analogs of (2.67)–(2.69), implying that the sequence $\{Y_{n\tau}\}$ is bounded and equicontinues on any bounded set of $\mathbb{R}^2$. We will prove now that any uniformly convergent subsequence $\{Y_{n_l\tau}\}$ has the same limit $Y_\tau$. To derive the limiting equation for $Y_\tau$, we treat $Y_{n\tau}$ as $Y_n$ of (2.63), and applying first the Duhamel formula (2.14), write

$$Y_{n\tau}(x,t) = \frac{i}{\sqrt{n}} \int_0^t \sum_{j,k=1}^n \mathbf{E}\{W_{jk}^{(n)\tau} \Phi_n\} \, dt_1$$

with $\Phi_n$ of (3.97) [cf. (2.70)]. Then an argument, similar to that leading to (3.63)–(3.67) and based on (3.6) with $p = 3$, yields

$$(3.101) \qquad Y_{n\tau}(x,t) = i \int_0^t \left( \sum_{l=0}^3 (T_l + r_l) + \varepsilon_{3\tau,n} \right) dt_1,$$

where [cf. (3.64)]

$$(3.102) \qquad T_l = \frac{1}{l! n^{(l+1)/2}} \sum_{j,k=1}^n \kappa_{l+1,jk} \mathbf{E}\{D_{jk}^l \Phi_n\}, \qquad l = 0,1,2,3,$$

$$(3.103) \qquad |\varepsilon_{3\tau,n}| \leq \tau C_4(t,x)$$

and $r_l$ satisfies (3.67) with $s = 1$. Besides, $T_0 = 0$, $T_2$ satisfies (3.55), and the contribution to $T_3$ due to the term $9w^4\delta_{jk}$ of $\kappa_{4,jk} = \kappa_4 - 9w^4\delta_{jk}$ [see (3.2) and (3.4)] is bounded by $C(t,x)n^{-1}$. This allows us to write an analog of (2.71) with additional term proportional to $\kappa_4$:

$$(3.104) \quad Y_{n\tau}(x,t) = T_{w^2,n}^\tau + T_{\kappa_4,n}^\tau + \mathcal{E}_{3\tau,n}(x,t) + o(1), \qquad n \to \infty,$$

where

$$(3.105) \qquad T_{w^2,n}^\tau = iw^2 \int_0^t \frac{1}{n} \sum_{j,k=1}^n (1 + \delta_{jk}) \mathbf{E}\{D_{jk} \Phi_n\} \, dt_1,$$

$$(3.106) \qquad T_{\kappa_4,n}^\tau = i\kappa_4 \int_0^t \frac{1}{6n^2} \sum_{j,k=1}^n \mathbf{E}\{D_{jk}^3 \Phi_n\} \, dt_1,$$

$$(3.107) \qquad \mathcal{E}_{3\tau,n}(t,x) = \int_0^t \varepsilon_{3\tau,n}(t_1,x) \, dt_1$$

and for any $\tau > 0$ the reminder term $o(1)$ in (3.104) vanishes as $n \to \infty$ uniformly on any compact set of $\{t \geq 0, x \in \mathbb{R}\}$. The term $T_{w^2,n}^\tau$ of (3.105) has the same form as the r.h.s. of (2.71) of the GOE case. Since the argument, leading from (2.71) to (2.76)–(2.78) does not use the Gaussian form of $W_{jk}$ in (2.71), it is applicable in our case as well and yields

$$T_{w^2,n}^\tau = -2w^2 \int_0^t dt_1 \int_0^{t_1} \overline{v}_{n\tau}(t_2) Y_{n\tau}(x, t_1 - t_2) \, dt_2$$
$$+ x Z_{n\tau}(x) A_{n\tau}(t) - r_{n\tau}(x, t),$$

where $\overline{v}_{n\tau} = n^{-1} \mathbf{E}\{u_{n\tau}\}$, and $A_{n\tau}$ and $r_{n\tau}$ are given by (2.77) and (2.78) with the GOE matrix $M$ replaced by the truncated Wigner matrix $M^\tau$. Now it follows from the Schwarz inequality, (3.69) and (3.40) that, for any $\tau > 0$, the reminder $r_{n\tau}(x, t)$ vanishes as $n \to \infty$ uniformly on any compact of $\{t \leq 0, x \in \mathbb{R}\}$. Besides, in view of $|\overline{v}_{n\tau}| \leq 1$ and (3.58), $\lim_{n \to \infty}(\overline{v}_{n\tau} - \overline{v}_n) = 0$, $\forall \tau > 0$, and Theorem 2.1 yields that, for any $\tau > 0$, the sequences $\{\overline{v}_{n\tau}\}$ and $\{A_{n\tau}\}$ converge uniformly as $n \to \infty$ on any finite interval of $\mathbb{R}$ to $v(t)$ and $A(t)$ of (2.80) and (2.81). It follows also from (3.93), Theorem 3.5 and (3.40) that

$$|Z'_{n\tau}(x)| \leq |x| \mathbf{Var}^{1/2}\{\mathcal{N}_{n\tau}[\varphi]\}$$
$$\leq |x| C_\tau(\mu_4) \int (1 + |t|^4) |\widehat{\varphi}(t)| \, dt < \infty.$$

Hence, the sequence $\{Z_{n\tau}\}_{n>0}$ is compact for any $\tau > 0$. Denoting the continuous limit of some its subsequence $\{Z_{n_l\tau}\}_{n>0}$ by $Z_\tau$, we have for any $\tau > 0$ uniformly on any compact set of $\{t \geq 0, x \in \mathbb{R}\}$

$$\lim_{n_l \to \infty} T_{w^2, n_l}^\tau = -2w^2 \int_0^t dt_1 \int_0^{t_1} v(t_2) Y_\tau(x, t_1 - t_2) \, dt_2 + x Z_\tau(x) A(t).$$

(3.108)

Consider now the term $T_{\kappa_4, n}^\tau$ of (3.106) and note first that, in view of (3.58) and (3.97), we can replace $T_{\kappa_4, n}^\tau$ by

$$(3.109) \qquad T_{\kappa_4, n} = i\kappa_4 \int_0^t \frac{1}{6n^2} \sum_{j,k=1}^n \mathbf{E}\{D_{jk}^3(U_{jk}(t_1) e_n^\circ(x))\} \, dt_1$$

with the error bounded by $C_4(t, x) \tau^{-4} L_n^{(4)}(\tau)$. It follows now from (2.17), (2.72), (3.54) and (3.55) that the contribution to $T_{\kappa_4, n}$ due to any term of

$$n^{-2} \sum_{j,k=1}^n D_{jk}^3(U_{jk}(t) e_n^\circ(x)),$$

containing at least one off-diagonal element $U_{jk}$, is bounded by $C_3(t, x) n^{-1}$. Thus, we are left with terms, containing only diagonal elements of $U$. These terms arise from $e_n^\circ(x) D_{jk}^3 U_{jk}(t)$ and $3 D_{jk} U_{jk}(t) D_{jk}^2 e_n^\circ(x)$ in the above sum, and, by (2.17) and (2.72), their contributions to $T_{\kappa_4, n}$ are

$$(3.110) \qquad \frac{\kappa_4}{n^2} \sum_{j,k=1}^n \int_0^t \mathbf{E}\{(U_{jj} * U_{jj} * U_{kk} * U_{kk})(t_1) e_n^\circ(x)\} \, dt_1$$

and

$$(3.111) \quad \frac{ix\kappa_4}{n^2} \sum_{j,k=1}^{n} \int_0^t dt_1 \int t_2 \widehat{\varphi}(t_2) \mathbf{E}\{(U_{jj} * U_{kk})(t_1)(U_{jj} * U_{kk})(t_2) e_n(x)\} \, dt_2,$$

where we omitted $\beta_{jk}^3$, because the corresponding error term is $O(n^{-1})$. It is easy to see that the entries of $U$ appear in (3.110) and (3.111) in the form

$$(3.112) \qquad\qquad \mathbf{E}\{v_n(t_1, t_2) v_n(t_3, t_4) e_n^\circ(x)\}$$

and

$$(3.113) \qquad\qquad \mathbf{E}\{v_n(t_1, t_2) v_n(t_3, t_4) e_n(x)\},$$

where

$$(3.114) \qquad\qquad v_n(t_1, t_2) = n^{-1} \sum_{j=1}^{n} U_{jj}(t_1) U_{jj}(t_2).$$

Since $|U_{jj}(t)| \leq 1, t \in \mathbb{R}$, we have

$$(3.115) \qquad\qquad\qquad |v_n(t_1, t_2)| \leq 1.$$

This, the inequality $|e_n^\circ(x)| \leq 2$, and the general inequality

$$(3.116) \qquad\qquad \mathbf{E}\{|(\xi_1\xi_2)^\circ|\} \leq 2c\mathbf{E}\{|\xi_1^\circ|\} + 2c\mathbf{E}\{|\xi_2^\circ|\},$$

where $\xi_{1,2}^\circ = \xi_{1,2} - \mathbf{E}\{\xi_{1,2}\}$, and $\xi_{1,2}$ are random variables such that $|\xi_{1,2}| \leq c$ allow us to write for (3.112)

$$|\mathbf{E}\{v_n(t_1, t_2) v_n(t_3, t_4) e_n^\circ(x)\}| \leq 4\mathbf{E}\{|v_n^\circ(t_1, t_2)|\} + 4\mathbf{E}\{|v_n^\circ(t_3, t_4)|\}.$$

By Lemma 3.1 below, we have

$$(3.117) \qquad\qquad \mathbf{E}\{|v_n^\circ(t_1, t_2)|\} \leq C_3(t_1, t_2) n^{-1/4}.$$

Thus, (3.112) vanishes as $n \to \infty$ uniformly in $t$ and $x$, varying in any compact set of $\mathbb{R}^2$.

Expression (3.113) can be written as the sum of (3.112) and

$$(3.118) \quad \mathbf{E}\{v_n(t_1, t_2) v_n(t_3, t_4)\} \mathbf{E}\{e_n(x)\} = \mathbf{E}\{v_n(t_1, t_2) v_n(t_3, t_4)\} Z_n(x).$$

It is follows from (3.115) and (3.117) that $\mathbf{E}\{v_n(t_1, t_2) v_n(t_3, t_4)\}$ can be written as the product $\overline{v}_n(t_1, t_2) \overline{v}_n(t_3, t_4)$ up to an error term bounded by $C_3(t_1, t_2)^{1/2} C_3(t_3, t_4)^{1/2} n^{-1/4}$, where

$$(3.119) \qquad\qquad \overline{v}_n(t_1, t_2) = \mathbf{E}\{v_n(t_1, t_2)\}.$$

In addition, we have, by Lemma 3.1 below,

$$(3.120) \qquad\qquad \overline{v}_n(t_1, t_2) = v(t_1) v(t_2) + o(1),$$

where $v$ is given by (2.80) and $o(1)$ is bounded by $C_3(t_1, t_2) n^{-1/2}$.

We conclude from the above that the contribution of (3.110) to $T_{\kappa_4, n_l}$ vanishes as $n_l \to \infty$ uniformly in $t$ and $x$, varying in any compact set of $\{t \geq 0, x \in \mathbb{R}\}$, while in (3.111) we can replace $U_{jj}$ and $U_{kk}$ by $v$. As a result, we obtain

$$
\begin{aligned}
(3.121) \quad T_{\kappa_4} &:= \lim_{n_l \to \infty} T_{\kappa_4, n_l} \\
&= ix Z_\tau(x) \kappa_4 \int_0^t (v * v)(t_1) \, dt_1 \int t_2 \widehat{\varphi}(t_2)(v * v)(t_2) \, dt_2,
\end{aligned}
$$

uniformly on any compact of $\{t \geq 0, x \in \mathbb{R}\}$.

In view of (2.83) and Proposition 2.1(iii), we have

$$
(v * v)(t) = -\frac{1}{2\pi} \int_L e^{itz} f^2(z) \, dz.
$$

The integral over $L$ can be replaced by that over the cut $[-2w, 2w]$ of $\sqrt{z^2 - 4w^2}$ in (2.36) and we obtain that

$$
(3.122) \qquad (v * v)(t) = -\frac{i}{2\pi w^4} \int_{-2w}^{2w} e^{it\mu} \mu \sqrt{4w^2 - \mu^2} \, d\mu
$$

or, integrating by parts,

$$
(3.123) \qquad \frac{1}{\pi t w^4} \int_{-2w}^{2w} e^{it\mu} \frac{2w^2 - \mu^2}{\sqrt{4w^2 - \mu^2}} \, d\mu.
$$

Now the Parseval equation implies that

$$
(3.124) \quad \int t \widehat{\varphi}(t)(v * v)(t) \, dt = \frac{1}{\pi w^4} \int_{-2w}^{2w} \varphi(\mu) \frac{2w^2 - \mu^2}{\sqrt{4w^2 - \mu^2}} \, d\mu =: B,
$$

thus,

$$
T_{\kappa_4} = iB I(t) x Z_\tau(x),
$$

where

$$
(3.125) \qquad I(t) = \int_0^t (v * v)(t_1) \, dt_1.
$$

Besides, it follows from (3.104) and the convergence of sequences $\{Y_{n_l \tau}\}$, $\{T^\tau_{w^2, n_l}\}$ and $\{T^\tau_{\kappa_4, n_l}\}$ [(3.108) and (3.121)] that the limit

$$
(3.126) \qquad \mathcal{E}_{3\tau}(t, x) = \lim_{n_l \to \infty} \mathcal{E}_{3\tau, n_l}(t, x)
$$

exists uniformly on any compact of $\{t \geq 0, x \in \mathbb{R}\}$, and we have by (3.103) and (3.107)

$$
(3.127) \qquad |\mathcal{E}_{3\tau}(t, x)| \leq \tau C_5(t, x).
$$

This and (3.108) allow us to pass to the limit $n_l \to \infty$ in (3.104) and to obtain the integral equation

$$(3.128) \quad \begin{aligned} & Y_\tau(x,t) + 2w^2 \int_0^t dt_1 \int_0^{t_1} v(t_1 - t_2) Y_\tau(x, t_2) \, dt_2 \\ & = x Z_\tau(x)[A(t) + i\kappa_4 BI(t)] + \mathcal{E}_{3\tau}(t,x). \end{aligned}$$

The l.h.s. of (3.128) coincides with that of (2.82) and the r.h.s. of (3.128) is equal to that of (2.82) plus two more terms. Thus, the solution of (3.128) is equal to the r.h.s. of (2.86) plus two more terms, the contributions of the additional terms in the r.h.s. of (3.128). The r.h.s. of (2.86) leads to the first term in (3.92) [see (2.87) and the subsequent argument]. To find the contribution to (3.92) of the second term of the r.h.s. of (3.128), we use the r.h.s. of (2.11) with $R(t) = i\kappa_4 x Z_\tau(x) BI(t)$, $T_1$ of (2.85), and (3.122). This leads to the term

$$\frac{ix Z_\tau(x) B}{2\pi w^4} \int_{-2w}^{2w} \frac{e^{it\lambda}(2w^2 - \lambda^2)}{\sqrt{4w^2 - \lambda^2}} \, d\lambda$$

in the solution of (3.128), where we used the relations

$$\int_{-2w}^{2w} \frac{1}{\sqrt{4w^2 - \lambda^2}(\lambda - \mu)} \, d\lambda = 0,$$

$$\int_{-2w}^{2w} \frac{\sqrt{4w^2 - \lambda^2}}{(\lambda - \mu)} \, d\lambda = -\pi\mu, \qquad |\mu| < 2w.$$

Then the limiting form of (2.62) and (2.57) yield the expression $\kappa_4 B^2/2$, that is, the second term of (3.92).

Let us consider the contribution $C_{3\tau}(t,x)$ of the third term of the r.h.s. of (3.128), which is given by the r.h.s. of (2.11) with $R(t) = \mathcal{E}_{3\tau}(t,x)$ and $T_1$ of (2.85). Integrating by parts, we obtain

$$(3.129) \quad C_{3\tau}(t,x) = T_1(0)\mathcal{E}_{3\tau}(t,x) + \int_0^t T_1'(t - t_1)\mathcal{E}_{3\tau}(t_1, x) \, dt_1.$$

We have also

$$T_1(t) = -J_0(2wt), \qquad T_1'(t) = 2w J_1(2wt),$$

where $J_0$ and $J_1$ are the corresponding Bessel functions, so that $|T_1(t)| \leq 1$, $|T_1'(t)| \leq 2w$, and

$$T_1'(t) = \sqrt{\frac{4w}{\pi t}} \sin(2wt - \pi/4)(1 + O(t^{-3/2})), \qquad t \to \infty.$$

By using this, (3.129) and (3.127), it can be shown that

$$(3.130) \qquad |C_{3\tau}(t,x)| \leq \tau C_5(t,x).$$

Now, the limiting form of (2.62) implies

$$\tag{3.131} Z'_\tau(x) = -x V_{\text{Wig}} Z_\tau(x) + D_{3\tau}(x),$$

where

$$D_{3\tau}(x) = i \int \widehat{\varphi}(t) C_{3\tau}(t, x) \, dt$$

and, in view of (3.130)

$$\tag{3.132} |D_{3\tau}(x)| \leq \tau C_4(x) \int (1 + |t|^5) |\widehat{\varphi}(t)| \, dt,$$

and $C_4$ is $n$- and $\tau$-independent polynomial in $|x|$ of degree 4.

Since $Z_\tau(0) = 1$, we can replace (3.131) by

$$Z_\tau(x) = e^{-V_{\text{Wig}} x^2 / 2} + \int_0^x e^{-V_{\text{Wig}} (x^2 - y^2)/2} D_{3\tau}(y) \, dy$$

and then (3.40) and (3.132) imply that

$$\lim_{\tau \to 0} Z_\tau(x) = e^{-V_{\text{Wig}} x^2 / 2},$$

hence the assertion of theorem. $\square$

REMARK 3.4. (1) The proof of the CLT under condition (3.89) is much simpler, because it does not use the truncation procedure.

(2) Another expression for the limiting variance of linear eigenvalue statistics is obtained in [2]. In fact, the paper deals with the more general class of random matrices that the authors called the band matrices and that includes the sample covariance matrices with uncorrelated entries of data matrices of Section 4 below. Thus, a rather general formula for the variance of linear eigenvalue statistics obtained in [2] reduces to formulas (4.28) and (4.65) below.

It remains to prove the following:

LEMMA 3.1. *We have under the assumptions of Theorem 3.6*

$$\tag{3.133} \mathbf{Var}\{v_n(t_1, t_2)\} \leq C_3(t_1, t_2)/n^{1/2}$$

*and*

$$\mathbf{E}\{v_n(t_1, t_2)\} = \overline{v}_n(t_1) \overline{v}_n(t_2) + r_n(t_1, t_2), \qquad |r_n(t_1, t_2)| \leq \frac{C_3(t_1, t_2)}{n^{1/2}},$$

$\tag{3.134}$

*where $C_p$ is a polynomial in $|t_j|$, $j = 1, 2$, of degree $p$ with positive coefficients.*

PROOF. We denote again $\widehat{M} = n^{-1/2}\widehat{W}$ the GOE matrix (2.23), $\widehat{U}(t) = e^{it\widehat{M}}$, $\widehat{v}_n(t_1, t_2) = n^{-1}\sum_{j=1}^{n}\widehat{U}_{jj}(t_1)\widehat{U}_{jj}(t_2)$, and write

$$\begin{aligned}
\mathbf{Var}\{v_n(t_1, t_2)\} &= \mathbf{E}\{\widehat{v}_n(t_1, t_2)v_n^{\circ}(-t_1, -t_2)\} \\
(3.135) &\qquad + \mathbf{E}\{(v_n(t_1, t_2) - \widehat{v}_n(t_1, t_2))v_n^{\circ}(-t_1, -t_2)\} \\
&=: R_1 + R_2.
\end{aligned}$$

The Poincaré inequality (2.17), (2.21) and (2.24) allow to obtain

$$(3.136) \qquad \mathbf{Var}\{\widehat{v}_n(t_1, t_2)\} \le 4w^2(t_1^2 + t_2^2)n^{-2}$$

and

$$(3.137) \qquad \mathbf{Var}\{\widehat{U}_{jj}(t)\} \le 2w^2t^2n^{-1},$$

hence,

$$(3.138) \qquad |R_1| \le C_1(t_1, t_2)n^{-1}.$$

To estimate $R_2$, we use again the interpolation matrix (3.20) and write

$$\begin{aligned}
R_2 &= \int_0^1 \frac{d}{ds}\frac{1}{n}\sum_{j=1}^{n}\mathbf{E}\{U_{jj}(s, t_1)U_{jj}(s, t_2)v_n^{\circ}(-t_1, -t_2)\}\,ds \\
&= \frac{i}{2n^{3/2}}\int_0^1\sum_{j,k=1}^{n}\mathbf{E}\{(s^{-1/2}W_{jk}^{(n)} - (1-s)^{-1/2}\widehat{W}_{jk})\Psi(s, t_1, t_2)\}\,ds,
\end{aligned}$$

where

$$\Psi(s, t_1, t_2) = v_n^{\circ}(-t_1, -t_2)(t_1 U_{jk}(s, t_1)U_{jj}(s, t_2) + t_2 U_{jk}(s, t_2)U_{jj}(s, t_1))$$

and $U(s, t)$ is defined in (3.46). We have by (2.20), (3.6) with $p = 1$, and (2.19)

$$|R_2| \le C_3(t_1, t_2)n^{-1/2}$$

[cf. (3.24)]. This, (3.138) and (3.135) yield (3.133).

To prove (3.134), we write

$$(3.139) \qquad \mathbf{E}\{v_n(t_1, t_2)\} = \mathbf{E}\{\widehat{v}_n(t_1, t_2)\} + \mathbf{E}\{v_n(t_1, t_2) - \widehat{v}_n(t_1, t_2)\},$$

where the second term similar to term $R_2$ above is modulo bounded by $C_3 \times (t_1, t_2)n^{-1/2}$.

It follows from the orthogonal invariance of the GOE that $\mathbf{E}\{\widehat{U}_{jj}(t)\} = \widehat{\overline{v}}_n(t)$, so that

$$(3.140) \quad \mathbf{E}\{\widehat{v}_n(t_1, t_2)\} = \widehat{\overline{v}}_n(t_1)\widehat{\overline{v}}_n(t_2) + \mathbf{E}\left\{n^{-1}\sum_{j=1}^{n}\widehat{U}_{jj}(t_1)\widehat{U}_{jj}^{\circ}(t_2)\right\},$$

where by (3.137) the second term is modulo bounded by $C_1(t_1, t_2)n^{-1}$. Besides, by using interpolation procedure, we can show that

$$(3.141) \qquad |\mathbf{E}\{\widehat{v}_n(t) - v_n(t)\}| \le C_3(t)n^{-1/2}.$$

This, (3.139) and (3.140) yield (3.134). $\quad\square$

## 4. Sample covariance matrices.

4.1. *Generalities.* We again confine ourselves to the real symmetric matrices. Thus, we consider in this section $n \times n$ real symmetric matrices

$$(4.1) \qquad M = Y^T Y, \qquad Y = n^{-1/2}X,$$

where $X = \{X_{\alpha j}^{(m,n)}\}_{\alpha,j=1}^{m,n}$ is the $m \times n$ real random matrix with the distribution

$$(4.2) \qquad \mathbf{P}_{mn}(dX) = \prod_{\alpha=1}^{m} \prod_{j=1}^{n} F_{\alpha j}^{(m,n)}(dX_{\alpha j}),$$

satisfying

$$\int X F_{\alpha j}^{(m,n)}(dX) = 0, \qquad \int X^2 F_{\alpha j}^{(m,n)}(dX) = a^2.$$

In other words, the entries $\{M_{jk}^{(m,n)}\}_{j,k=1}^{n}$ of $M$ of (4.1) are

$$M_{jk}^{(m,n)} = n^{-1} \sum_{\alpha=1}^{m} X_{\alpha j}^{(m,n)} X_{\alpha k}^{(m,n)},$$

where $X_{\alpha j}^{(m,n)} \in \mathbb{R}$, $\alpha = 1, \ldots, m$, $j = 1, \ldots, n$, are independent random variables such that

$$(4.3) \qquad \mathbf{E}\{X_{\alpha j}^{(m,n)}\} = 0, \qquad \mathbf{E}\{X_{\alpha j}^{(m,n)} X_{\beta k}^{(m,n)}\} = \delta_{\alpha\beta}\delta_{jk}a^2.$$

A particular case of (4.1)–(4.3),

$$(4.4) \qquad \widehat{M} = \widehat{Y}^T\widehat{Y}, \qquad \widehat{Y} = n^{-1/2}\widehat{X},$$

where the entries of $\widehat{X} = \{\widehat{X}_{\alpha j}\}_{\alpha,j=1}^{m,n}$ are i.i.d. Gaussian random variables satisfying (4.3), that is,

$$(4.5) \qquad \mathbf{P}(d\widehat{X}) = \widehat{Z}_{mn1}^{-1} \exp\{-\operatorname{Tr}\widehat{X}^T\widehat{X}/2a^2\} \prod_{\alpha=1}^{m} \prod_{j=1}^{n} d\widehat{X}_{\alpha j},$$

is closely related to the null (white) case of the Wishart random matrix of statistics (see [21], Section 3.2). The difference is in the factor $m^{-1/2}$ instead of $n^{-1/2}$ in (4.4). In what follows, to simplify the notation, we will often omit the superscript $(m, n)$, and the sums over the Latin indexes will be from 1 to $n$, and the sums over the Greek indexes will be from 1 to $m$.

We present first an analog of the law of large numbers for the sample covariance matrices.

THEOREM 4.1. *Let $M$ be the real symmetric sample covariance matrix* (4.1)–(4.3). *Assume that for any $\tau > 0$*

$$(4.6) \qquad \frac{1}{n^2} \sum_{\alpha, j} \int_{|X| > \tau \sqrt{n}} X^2 F_{\alpha j}^{(m,n)}(dX) \to 0$$

*as*

$$(4.7) \qquad n \to \infty, \ m \to \infty, \qquad m/n \to c \in [0, \infty),$$

*and that $\{X_{\alpha j}\}_{\alpha, j=1}^{m,n}$ are defined on the same probability space for all $m, n \in \mathbb{N}$. Then for any bounded continuous $\varphi : \mathbb{R} \to \mathbb{C}$, we have with probability 1*

$$(4.8) \qquad \lim_{m, n \to \infty, m/n \to c} n^{-1} \mathcal{N}_n[\varphi] = \int \varphi(\lambda) N_{MP}(d\lambda),$$

*where $\mathcal{N}_n[\varphi]$ is defined in* (1.1), *and*

$$(4.9) \quad N_{MP}(d\lambda) = (1 - c)_+ \delta_0(\lambda) \, d\lambda + (2\pi a^2 \lambda)^{-1} \sqrt{((\lambda - a_-)(a_+ - \lambda))_+} \, d\lambda$$

*with $a_\pm = a^2 (1 \pm \sqrt{c})^2$, and $x_+ = \max(x, 0)$.*

We refer the reader to [4] and [12] for results and references concerning this assertion that dates back to [20]. Here we outline a weaker version of the theorem on the convergence in mean in (4.8), basing it on the same ideas as in Theorems 2.1–3.2. We will need this assertion as well as the method of its proof.

We start from the Gaussian case, that is, the Wishart random matrices (4.4) and (4.5), and follow essentially the proof of Theorem 2.1. Introduce [cf. (2.37)]

$$(4.10) \quad f_n(z) = \mathbf{E}\{g_n(z)\}, \qquad g_n(z) = n^{-1} \operatorname{Tr} G(z), \qquad G(z) = (M - z)^{-1}.$$

By using again the resolvent identity (2.40) and the Gaussian differentiation formula (2.20), we obtain for $f_n$ of (2.37)

$$f_n(z) = -\frac{1}{z} + \frac{1}{z n^{3/2}} \sum_{\alpha, k} \mathbf{E}\{\widehat{X}_{\alpha k} (\widehat{Y} G)_{\alpha k}(z)\}$$

$$= -\frac{1}{z} + \frac{a^2}{z n^2} \sum_{\alpha, k} \mathbf{E}\{\widehat{D}_{\alpha k} (\widehat{Y} G)_{\alpha k}(z)\}, \qquad \widehat{D}_{\alpha k} = \partial / \partial \widehat{Y}_{\alpha k}.$$

We have from (2.40) and (4.4) [cf. (2.42)]

$$(4.11) \qquad \widehat{D}_{\alpha k} G_{jk} = -(\widehat{Y}G)_{\alpha k} G_{jk} - (\widehat{Y}G)_{\alpha j} G_{kk},$$

hence, by (2.44) [cf. (2.43)],

$$(4.12) \qquad \begin{aligned} f_n(z) &= -\frac{1}{z} + \frac{ma^2}{nz} f_n(z) - \frac{a^2}{z} \mathbf{E}\{g_n(z)n^{-1}\operatorname{Tr} G(z)\widehat{M}\} \\ &\quad - \frac{a^2}{z}\mathbf{E}\{n^{-2}\operatorname{Tr} G^2(z)\widehat{M}\} \end{aligned}$$

or, after using the identity

$$(4.13) \qquad G(z)\widehat{M} = zG(z) + 1,$$

we get

$$\begin{aligned} f_n(z) &= -\frac{1}{z} + \frac{a^2}{z} c_n f_n(z) - \frac{a^2}{z}\mathbf{E}\{g_n(z)(zg_n(z)+1)\} \\ &\quad - \frac{a^2}{nz}\mathbf{E}\{n^{-1}\operatorname{Tr} G(z)(zG(z)+1)\} \end{aligned}$$

with

$$(4.14) \qquad c_n = m/n.$$

We need now the Poincaré type inequalities for the Wishart matrices (4.4) and (4.5):

$$(4.15) \qquad \mathbf{Var}\{\mathcal{N}_n[\varphi]\} \leq 4a^2 \mathbf{E}\{n^{-1}\operatorname{Tr} \varphi'(\widehat{M})\overline{\varphi'(\widehat{M})}\widehat{M}\}$$

$$(4.16) \qquad \leq 4a^4 c_n \sup_{\lambda \in \mathbb{R}} |\varphi'(\lambda)|^2$$

[cf. (2.25) and (2.26)]. They can be easily derived from (2.21) by using the formulas $\mathbf{E}\{n^{-1}\operatorname{Tr}\widehat{M}\} = a^2 c_n$ [see (4.3)] and

$$(4.17) \qquad \widehat{D}_{\alpha k} U_{jl}(t) = i(((\widehat{Y}U)_{\alpha j} * U_{kl})(t) + ((\widehat{Y}U)_{\alpha l} * U_{jk})(t)),$$

$$(4.18) \qquad \widehat{D}_{\alpha k} \operatorname{Tr} \varphi(\widehat{M}) = 2(\widehat{Y}\varphi'(\widehat{M}))_{\alpha k}$$

[cf. (2.17) and (2.28)]. By applying (4.16) with $\varphi(\lambda) = n^{-1}(\lambda - z)^{-1}$, we obtain [cf. (2.45)]

$$(4.19) \qquad \mathbf{Var}\{g_n(z)\} \leq \frac{4a^4 c_n}{n^2 |\Im z|^4}.$$

Thus, (2.40) and (4.12) allow us to write

$$f_n(z) = -z^{-1} - a^2 z^{-1}(1-c_n)f_n(z) - a^2 f_n^2(z) + r_n,$$

where

$$|r_n| \leq \frac{4a^6 c_n}{n^2 |\Im z|^5} + \frac{2a^2}{n |\Im z|^2}.$$

This and the limit (4.7) yield an analog of (2.46):

$$(4.20) \qquad za^2 f_{MP}^2(z) + (z + a^2(1-c))f_{MP}(z) + 1 = 0, \qquad \Im z \neq 0,$$

and since $\Im f_{MP}(z) \Im z \geq 0$, we obtain [cf. (2.36)]

$$(4.21) \qquad f_{MP}(z) = (\sqrt{(z - a_m)^2 - 4a^4 c} - (z + a^2(1-c)))/2a^2 z,$$

where the branch of the square root is fixed by the asymptotic form $z + O(1)$, $z \to \infty$, and

$$(4.22) \qquad\qquad\qquad a_m = a^2(c+1).$$

This and inversion formula

$$N_{MP}(\Delta) = \lim_{\varepsilon \to 0} \pi^{-1} \int_\Delta \Im f_{MP}(\lambda + i0) \, d\lambda,$$

where the endpoints of $\Delta$ are not the atoms of $N_{MP}$, lead to (4.9).

The next step is to prove an analog of Theorem 3.1, assuming that

$$(4.23) \qquad a_3 := \sup_n \max_{1 \leq \alpha \leq m, 1 \leq j \leq n} \mathbf{E}\{|X_{\alpha j}^{(m,n)}|^3\} < \infty.$$

To this end, we use again an "interpolation" matrix [cf. (3.20)]

$$(4.24) \qquad M(s) = Y^T(s)Y(s), \qquad Y(s) = s^{1/2}Y + (1-s)^{1/2}\widehat{Y}, \qquad s \in [0, 1],$$

where $Y$ and $\widehat{Y}$ are defined in (4.1)–(4.5). We have, with the same notation as in Theorem 3.1,

$$f_n(z) - \widehat{f}_n(z)$$

$$(4.25) \qquad = \int_0^1 \frac{\partial}{\partial s} \mathbf{E}\{n^{-1} \operatorname{Tr} G(s, z)\} \, ds$$

$$= -n^{-3/2} \int_0^1 \sum_{\alpha, k} \mathbf{E}\{(s^{-1/2} X_{\alpha k} - (1-s)^{-1/2} \widehat{X}_{\alpha k})(Y(s)G')_{\alpha k}\} \, ds.$$

Since $\{X_{\alpha j}\}_{\alpha,j=1}^{m,n}$ are independent random variables satisfying (4.3) and (4.23), and $\{\widehat{X}_{\alpha j}\}_{\alpha,j=1}^{m,n}$ are i.i.d. Gaussian random variables also satisfying (4.3), we apply the general differentiation formula (3.6) with $\Phi = (Y(s)G')_{\alpha j}$ and $p = 1$ to the contribution of the first term in the parentheses of (4.25) and the Gaussian differentiation formula to the contribution of the second term. As it was already several times in the case of the Wigner matrices (see, e.g., Corollary 3.1 and Theorem 3.1), the term with the first derivative of the general differentiation formula is canceled by the expression resulting from the Gaussian differentiation formula, and we are left with [cf. (3.24)]

$$f_n(z) - \widehat{f}_n(z) = \int_0^1 \sqrt{s}\,\varepsilon_1(s)\,ds,$$

where

$$|\varepsilon_1(s)| \leq \frac{C_1 a_3}{n^{5/2}} \sum_{\alpha,k} \sup_{Y \in \mathcal{M}_{m,n}} |D_{\alpha j}^2(YG')_{\alpha k}|, \qquad D_{\alpha k} = \partial/\partial Y_{\alpha k},$$

where $G = (Y^T Y - z)^{-1}$, and $\mathcal{M}_{m,n}$ is the set of $m \times n$ real matrices. It suffices to find an $O(1)$ bound for $D_{\alpha k}^2(YG')_{\alpha k}$. Since $(YG)_{\alpha k}$ is analytic in $z$, $\Im z \neq 0$, then the bound for $(YG')_{\alpha k}$ follows from that for $(YG)_{\alpha k}$ and the Cauchy bound for derivatives of analytic function. By using (4.11) and a little algebra, we obtain

$$D_{\alpha k}^2(YG)_{\alpha k} = -6G_{kk}(YG)_{\alpha k} + 6G_{kk}(YG)_{\alpha k}(YGY^T)_{\alpha\alpha} + 2(YG)_{\alpha k}^3.$$

It follows from (2.40) that $|G_{kk}| \leq |\Im z|^{-1}$. Next, if $G = (Y^T Y - z)^{-1}$ and $\widetilde{G} = (YY^T - z)^{-1}$, then $YG = \widetilde{G}Y$, and $(YGY^T)_{\alpha\alpha} = (\widetilde{G}YY^T)_{\alpha\alpha} = (1 + z\widetilde{G})_{\alpha\alpha}$ [see (4.13)], thus,

$$(4.26) \qquad |(YGY^T)_{\alpha\alpha}| \leq 1 + |z||\Im z|^{-1}.$$

Furthermore, it follows from the Schwarz inequality that

$$(4.27) \qquad |(YG)_{\alpha k}| \leq (G^* Y^T YG)_{kk}^{1/2} \leq ((1 + |z||\Im z|^{-1})/|\Im z|)^{1/2}.$$

Thus, $D_{\alpha k}^2(YG)_{\alpha k}$ is bounded uniformly in $1 \leq \alpha \leq m$, $1 \leq k \leq n$, all $m$ and $n$, and $z$, varying in a compact set $K \subset \mathbb{C} \setminus \mathbb{R}$, and

$$|\varepsilon_1(s)| \leq C_K n^{-1/2}, \qquad n \to \infty, m \to \infty, m/n \to c \in [0, \infty),$$

where $C_K < \infty$ depends only on $K \subset \mathbb{C} \setminus \mathbb{R}$.

In fact, a bit more tedious algebra and (4.26) and (4.27) show that, for every $1 \leq \alpha \leq m$, $1 \leq k \leq n$ $(YG)_{\alpha k}(z)$ is real analytic in every $Y_{\beta j}$, $1 \leq \beta \leq m$, $1 \leq j \leq n$ and $\Im z \neq 0$. Hence, all derivatives $\partial_{\alpha k}^l(YG)_{\alpha k}(z)$, $l = 0, 1, \ldots$, are bounded by $C_{lK}$, $z \in K \subset \mathbb{C} \setminus \mathbb{R}$ [cf. (3.25)].

This proves (4.9) under condition (4.23), that is, an analog of Theorem 3.1. To prove (4.9) under condition (4.6), we have to use the truncation procedure analogous to that of the proof of Theorem (3.2) and bounds (2.40), (4.26) and (4.27).

4.2. *Central limit theorem for linear eigenvalue statistics of the Wishart ensemble.* Following our scheme of the presentation in the case of the Wigner matrices, we start from the central limit theorem for the sample covariance matrices with Gaussian entries, that is, from the Wishart ensemble (4.4) and (4.5). We confine ourselves to the case $c \geq 1$.

THEOREM 4.2. *Let $\mathcal{N}_n[\varphi]$ be a linear eigenvalue statistic of the Wishart matrix (4.4) and (4.5), corresponding to a bounded function $\varphi : \mathbb{R} \to \mathbb{R}$ with bounded derivative. Then the centered random variable $\mathcal{N}_n^\circ[\varphi]$ [see (2.47)] converges in distribution as $n \to \infty, m \to \infty, m/n \to c \geq 1$ to the Gaussian random variable of zero mean and variance*

$$V_{\mathrm{Wish}}[\varphi] = \frac{1}{2\pi^2} \int_{a_-}^{a_+} \int_{a_-}^{a_+} \left( \frac{\triangle \varphi}{\triangle \lambda} \right)^2$$

$$(4.28) \qquad \times \frac{4a^4 c - (\lambda_1 - a_m)(\lambda_2 - a_m)}{\sqrt{4a^4 c - (\lambda_1 - a_m)^2}\sqrt{4a^4 c - (\lambda_2 - a_m)^2}} \, d\lambda_1 \, d\lambda_2,$$

*where $\Delta \varphi$ is defined in (2.49), $a_\pm = a^2(1 \pm \sqrt{c})^2$ and $a_m$ is defined in (4.22).*

PROOF. We follow the scheme of the proof of Theorem 2.2. Namely, assume first that $\varphi$ admits the Fourier transform $\widehat{\varphi}$ [see (2.53)], satisfying (2.54). We have, similarly to the proof of Theorem 2.2, the relations (2.57)–(2.64). It follows also from (4.16) with $\varphi(\lambda) = e^{it\lambda}$ that [cf. (2.65)]

$$(4.29) \qquad\qquad \mathbf{Var}\{u_n(t)\} \leq 4a^4 t^2 c_n,$$

thus [cf. (2.67)],

$$(4.30) \qquad |Y_n(x,t)| = |\mathbf{E}\{u_n^\circ(t)e_n(x)\}| \leq \mathbf{Var}^{1/2}\{u_n(t)\} \leq 2a^2|t|c_n^{1/2}.$$

Likewise, we have the bound

$$(4.31) \qquad\qquad |\partial Y_n(x,t)/\partial x| \leq 4a^4\sqrt{c_n}\sup_{\lambda \in \mathbb{R}}|\varphi'(\lambda)|,$$

following from (2.21) and (4.29) [cf. (2.69)], and the bound

$$(4.32) \qquad\qquad |\partial Y_n(x,t)/\partial t| \leq 2a^2\sqrt{c_n}(1 + Ca^4t^2)^{1/2}$$

with $C$ depending only on $c_n$, following from (4.5) and (4.15) [cf. (2.68)]. Hence, the sequence $\{Y_n\}$ is bounded and equicontinuous on any finite set of $\mathbb{R}^2$. We will prove now that any uniformly converging subsequence of $\{Y_n\}$ has the same limit $Y$, leading to (2.55), hence to (2.51) and (2.52) with $V_{\mathrm{Wish}}[\varphi]$ instead of $V_{\mathrm{GOE}}[\varphi]$. Applying the Duhamel formula (2.14), (2.20),

(4.17) and (4.18), we obtain

$$Y_n(x,t) = ia^2 c_n \int_0^t Y_n(x,t_1)\,dt_1 - a^2 n^{-1} \int_0^t dt_1 \int_0^{t_1} \mathbf{E}\{\operatorname{Tr}\widehat{M}U(t_1)e_n^\circ(x)\}\,dt_2$$

$$- a^2 n^{-1} \int_0^t dt_1 \int_0^{t_1} \mathbf{E}\{\operatorname{Tr}\widehat{M}U(t_1-t_2)\operatorname{Tr}U(t_2)e_n^\circ(x)\}\,dt_2$$

$$- 2a^2 x n^{-1} \int_0^t \mathbf{E}\{\operatorname{Tr}\varphi'(\widehat{M})\widehat{M}U(t_1)e_n(x)\}\,dt_1$$

or

$$\begin{aligned}
Y_n(x,t) = {}& ia^2(c_n-1)\int_0^t Y_n(x,t_1)\,dt_1 \\
& + ia^2 n^{-1}\int_0^t \mathbf{E}\{u_n'(t_1)e_n^\circ(x)\}t_1\,dt_1 \\
& + ia^2 n^{-1}\int_0^t \mathbf{E}\{u_n(t-t_1)u_n(t_1)e_n^\circ(x)\}\,dt_1 \\
& + 2ia^2 x n^{-1}\mathbf{E}\{\operatorname{Tr}\varphi'(\widehat{M})(U(t)-1)e_n(x)\},
\end{aligned} \tag{4.33}$$

where we used the formulas $\operatorname{Tr}\widehat{M}U(t) = -iu_n'(t)$ and

$$\int_0^t dt_1 \int_0^{t_1}\mathbf{E}\{u_n'(t_1-t_2)u_n(t_2)e_n^\circ(x)\}\,dt_2$$

$$= \int_0^t \mathbf{E}\{(u_n(t-t_1)-n)u_n(t_1)e_n^\circ(x)\}\,dt_1.$$

This and an analog of (2.74) and (2.75) yield an analog of (2.76),

$$Y_n(x,t) - ia^2(c_n-1)\int_0^t Y_n(x,t_1)\,dt_1 - 2ia^2\int_0^t \overline{v}_n(t-t_1)Y_n(x,t_1)\,dt_1$$

$$= 2ia^2 x Z_n(x)\int \varphi'(\lambda)(e^{it\lambda}-1)\mathbf{E}\{N_n(d\lambda)\} + r_n(x,t), \tag{4.34}$$

where now $\overline{v}_n = n^{-1}\mathbf{E}\{u_n\}$ and

$$\begin{aligned}
r_n(x,t) = {}& ia^2 n^{-1}\int_0^t (Y_n(x,t)-Y_n(x,t_1))\,dt_1 \\
& + ia^2 n^{-1}\int_0^t \mathbf{E}\{u_n^\circ(t-t_1)u_n^\circ(t_1)e_n^\circ(x)\}\,dt_1 \\
& - 2a^2 x n^{-1}\int \theta\widehat{\varphi}(\theta)(Y_n(x,t+\theta)-Y_n(x,\theta))\,d\theta. 
\end{aligned} \tag{4.35}$$

It follows from (2.54), (4.29) and (4.30) that $r_n(x,t) = O(n^{-1})$ uniformly in $(x,t)$, varying in a compact set $K \subset \{x \in \mathbb{R}, t \geq 0\}$. This and Theorem 4.1 imply that the limit of every uniformly converging subsequence of $\{Y_n\}$ solves the equation [cf. (2.82)]

$$Y(x,t) - ia^2(c-1)\int_0^t Y(x,t_1)\,dt_1$$

$$-2ia^2\int_0^t v_{MP}(t-t_1)Y(x,t_1)\,dt_1 = xZ(x)A(t),$$

where [cf. (2.80)]

$$(4.36) \quad v_{MP}(t) := \lim_{n\to\infty}\overline{v}_n(t) = \frac{1}{2\pi a^2}\int_{a_-}^{a_+} e^{it\lambda}\sqrt{4a^2c-(\lambda-a_m)^2}\,\lambda^{-1}\,d\lambda,$$

and [cf. (2.81)]

$$(4.37) \quad \begin{aligned} A(t) &= 2a^2i\int \varphi'(\lambda)(e^{it\lambda}-1)N_{MP}(d\lambda) \\ &= -\frac{1}{\pi}\int_0^t dt_1\int_{a_-}^{a_+} e^{it_1\lambda}\varphi'(\lambda)\sqrt{4a^4c-(\lambda-a_m)^2}\,d\lambda. \end{aligned}$$

Now an argument similar to that leading from (2.82) to (2.86) and based on Proposition 2.1 and the formula

$$(4.38) \qquad\qquad \widehat{v}_{MP} = f_{MP}$$

yields

$$(4.39) \quad \begin{aligned} Y(x,t) &= \frac{ixZ(x)}{\pi^2}\int_{a_-}^{a_+}\varphi'(\lambda)\sqrt{4a^4c-(\mu-a_m)^2}\,d\lambda \\ &\quad \times \int_{a_-}^{a_+}\frac{e^{it\mu}-e^{it\lambda}}{\sqrt{4a^4c-(\mu-a_m)^2}(\mu-\lambda)}\,d\mu. \end{aligned}$$

Using this in (2.62), we obtain an analog of (2.87), and then an analog of (2.55) via (2.58) with $V_{\text{Wish}}$ of (4.28) instead of $V_{\text{GOE}}$, that is, an equation for the limiting characteristic function. Since the equation is uniquely soluble, we have finally

$$Z(x) = e^{-x^2 V_{\text{Wish}}[\varphi]/2},$$

that is, the assertion of the theorem under condition (2.54). The general case of bounded test functions with bounded derivative can be obtained via an approximation procedure analogous to that of the end of the proof of Theorem 2.2 and based on (4.16).  $\square$

REMARK 4.1.   (1) The proof of Theorem 4.2 can be easily modified to obtain an analogous assertion for the Laguerre ensemble of Hermitian matrices $M = n^{-1}X^*X$, where the complex $m \times n$ matrix $X$ has the probability

distribution [cf. (4.5)]

$$\mathbf{P}(dX) = Z_{mn2}^{-1} \exp\{-\operatorname{Tr} X^* X / a^2\} \prod_{\alpha=1}^{m} \prod_{j=1}^{n} d\Re X_{\alpha j} \, d\Im X_{\alpha j}.$$

The result is given by Theorem 4.2, in which $V_{\mathrm{Wish}}$ is replaced by $V_{\mathrm{Lag}} = V_{\mathrm{Wish}}/2$.

(2) It follows from the representation of the density $\rho_n$ of $\mathbf{E}\{N_n\}$ via the Laguerre polynomials that (see [19], Chapters 6 and 7)

$$\rho_n(\lambda) \leq C e^{-cn\lambda}$$

for finite $c$ and $C$ and $\lambda$ sufficiently big. This bound and the approximation procedure of the end of proof of Theorem 2.2 allows us to extend the theorem to $C^1$ test functions whose derivative grows as $C_1 e^{c_1 \lambda}$ for some $c_1 > 0$ and $C_1 < \infty$.

### 4.3. Central limit theorem for linear eigenvalue statistics of sample covariance matrices: the case of zero excess of entries. We prove here an analog of Theorem 3.4 for the sample covariance matrices.

THEOREM 4.3. *Let $M$ be the sample covariance matrix (4.1)–(4.3). Assume the following:*

(i) *the third and fourth moments of entries do not depend on $j$, $k$, $m$ and $n$:*

$$(4.40) \qquad \mu_3 = \mathbf{E}\{(X_{\alpha j}^{(m,n)})^3\}, \qquad \mu_4 = \mathbf{E}\{(X_{\alpha j}^{(m,n)})^4\};$$

(ii) *for any $\tau > 0$,*

$$(4.41) \qquad L_{mn}^{(4)}(\tau) := n^{-2} \sum_{\alpha,j} \int_{|X|>\tau\sqrt{n}} X^4 F_{\alpha j}^{(m,n)}(dX) \to 0$$

*as $n \to \infty$, $m \to \infty$, $m/n \to c \in [1,\infty)$;*

(iii) *the fourth cumulant of entries is zero:*

$$(4.42) \qquad \kappa_4 = \mu_4 - 3a^4 = 0.$$

*Let $\varphi : \mathbb{R} \to \mathbb{R}$ be a test function whose Fourier transform satisfies (3.40).*

*Then the corresponding centered linear eigenvalue statistic $\mathcal{N}_n^{\circ}[\varphi]$ converges in distribution to the Gaussian random variable of zero mean and variance $V_{\mathrm{Wish}}[\varphi]$ of (4.28).*

PROOF. We follow the scheme of the proof of Theorem 3.4. Thus, in view of Theorem 4.2, it suffices to prove that if subsequently

$$(4.43) \qquad m, n \to \infty, \qquad m/n \to c \in [1,\infty) \quad \text{and} \quad \tau \to 0,$$

then [cf. (3.60)]

$$(4.44) \qquad R^\tau_{mn}(x) = \mathbf{E}\{e^{ix\mathcal{N}^\circ_{n\tau}[\varphi]}\} - \mathbf{E}\{e^{ix\widehat{\mathcal{N}}^\circ_n[\varphi]}\} \to 0,$$

where $\mathcal{N}_{n\tau}[\varphi]$ is a linear eigenvalue statistic corresponding to the truncated matrix [cf. (3.27)]

$$(4.45) \qquad M^\tau = (Y^\tau)^T Y^\tau, \qquad Y^\tau = n^{-1/2} X^\tau,$$

$$X^\tau = \{X^\tau_{\alpha j} = \text{sign}\, X^{(m,n)}_{\alpha j} \max\{|X^{(m,n)}_{\alpha j}|, \tau n^{1/2}\}\}^{m,n}_{\alpha, j=1}$$

and the statistic $\widehat{\mathcal{N}}_n[\varphi]$ corresponds to the Wishart matrix $\widehat{Y}^T \widehat{Y}$ of (4.4). By using interpolating matrix [cf. (3.35) and (4.24)]

$$(4.46) \qquad M^\tau(s) = Y^{\tau T}(s) Y^\tau(s),$$

$$Y^\tau(s) = s^{1/2} Y^\tau + (1-s)^{1/2}\widehat{Y}, \qquad s \in [0,1],$$

we have [cf. (3.43)–(3.46)]

$$(4.47) \qquad R^\tau_{mn}(x) = -x \int_0^1 ds \int t\widehat{\varphi}(t)[A_n - B_n]\, dt,$$

where now

$$A_n = \frac{1}{\sqrt{ns}} \sum_{\alpha,k} \mathbf{E}\{X^\tau_{\alpha k} \Phi_{\alpha k}(s)\}, \qquad B_n = \frac{1}{\sqrt{n(1-s)}} \sum_{\alpha,k} \mathbf{E}\{\widehat{X}_{\alpha k} \Phi_{\alpha k}(s)\}$$

with

$$(4.48) \qquad \Phi_{\alpha k}(s) = e^\circ_n(s,x)(Y^\tau(s)U(s,t))_{\alpha k}$$

and $e_n(s,x)$ and $U(s,t)$ are defined in (3.62) in which $M^\tau$ is given by (4.45) and (4.46). We have, by (2.20),

$$B_n = \frac{a^2}{n} \sum_{\alpha,k} \mathbf{E}\{D_{\alpha k}(s)\Phi_{\alpha k}(s)\}, \qquad D_{\alpha k}(s) = \partial/\partial Y^\tau_{\alpha k}(s)$$

and, by (3.6) with $p = 3$ [cf. (3.63)–(3.65)],

$$A_n = \sum_{l=0}^3 T_{l\tau} + \varepsilon_{3\tau},$$

where now

$$(4.49) \quad T_{l\tau} = \frac{s^{(l-1)/2}}{l!n^{(l+1)/2}} \sum_{\alpha,k} \kappa^\tau_{l+1,\alpha k} \mathbf{E}\{D^l_{\alpha k}(s)\Phi_{\alpha k}(s)\}, \qquad l = 0,1,2,3,$$

$\kappa^\tau_{l,\alpha k}$ is $l$th cumulant of $X^\tau_{\alpha k}$, and

$$|\varepsilon_{3\tau}| \le \frac{C_3\mu_4\tau}{n^2} \sum_{\alpha,k} \sup_{|X| \le \tau\sqrt{n}} |\mathbf{E}\{D^4_{\alpha k}(s)\Phi_{\alpha k}(s)|_{Y^\tau_{\alpha k}(s)=(s/n)^{1/2}X+(1-s)^{1/2}\widehat{Y}_{\alpha k}}\}|$$

in view of $\mathbf{E}\{|X^{\tau}_{\alpha k}|^5\} \leq \tau \sqrt{n} \mu_4$.

In what follows we omit $s$ and denote $D_{\alpha k} = D_{\alpha k}(s)$, $U(t) = U(t, s)$, etc. Let us prove the uniform boundedness of derivatives $\mathbf{E}\{D^l_{\alpha k}\Phi_{\alpha k}\}$, $l \leq 4$, that will allow us to obtain analogs of (3.65)–(3.67). We have, analogously to (4.17),

$$(4.50) \qquad D_{\alpha k}U_{jk}(t) = i(((Y^{\tau}U)_{\alpha k} * U_{jk})(t) + ((Y^{\tau}U)_{\alpha j} * U_{kk})(t)),$$

$$(4.51) \qquad D_{\alpha k}e_n(x) = -2x e_n(x) \int \theta \widehat{\varphi}(\theta)(Y^{\tau}U)_{\alpha k}(\theta)\, d\theta,$$

$$(4.52) \qquad \begin{aligned} D_{\alpha k}(Y^{\tau}U)_{\alpha k}(t) = U_{kk}(t) &+ i(((Y^{\tau}UY^{\tau T})_{\alpha\alpha} * U_{kk})(t) \\ &+ ((Y^{\tau}U)_{\alpha k} * (Y^{\tau}U)_{\alpha k})(t)), \end{aligned}$$

$$(4.53) \quad D_{\alpha k}(Y^{\tau}UY^{\tau T})_{\alpha\alpha}(t) = 2(Y^{\tau}U)_{\alpha k}(t) + 2i((Y^{\tau}UY^{\tau T})_{\alpha\alpha} * (Y^{\tau}U)_{\alpha k})(t).$$

Since by (2.16)

$$(4.54) \quad |(Y^{\tau}U)_{\alpha k}(t)| \leq \left(\sum_j (Y^{\tau}_{\alpha j})^2\right)^{1/2}, \qquad |(Y^{\tau}UY^{\tau T})_{\alpha\alpha}(t)| \leq \sum_j (Y^{\tau}_{\alpha j})^2,$$

then, iterating (4.50)–(4.53), we have

$$(4.55) \qquad |D^l_{\alpha k}\Phi_{\alpha k}| \leq C_l(t, x)\left(\sum_j (Y^{\tau}_{\alpha j})^2\right)^{(l+1)/2}$$

and by (4.46),

$$|\mathbf{E}\{D^l_{\alpha k}\Phi_{\alpha k}\}| \leq \frac{C_l(t, x)}{n^{(l+1)/2}}\mathbf{E}\left\{\left(\sum_j (X^{\tau}_{\alpha j})^2\right)^{(l+1)/2} + \left(\sum_j (\widehat{X}_{\alpha j})^2\right)^{(l+1)/2}\right\},$$
$$l \geq 0.$$

Now the Hölder inequality implies the bound

$$n^{-(l+1)/2}\mathbf{E}\left\{\left(\sum_j (X^{\tau}_{\alpha j})^2\right)^{(l+1)/2}\right\} \leq n^{-1}\mathbf{E}\left\{\sum_j |X^{\tau}_{\alpha j}|^{l+1}\right\} \leq \mu_4^{(l+1)/4}, \qquad l \leq 3,$$

and analogous bounds for $\{\widehat{X}_{\alpha j}\}$, thus,

$$(4.56) \qquad |\mathbf{E}\{D^l_{\alpha k}\Phi_{\alpha k}\}| \leq C_l(t), \qquad l \leq 3.$$

In the case where $l = 4$ a similar argument and (4.55) yield

$$\sup_{|X| \leq \tau\sqrt{n}} |\mathbf{E}\{D^4_{\alpha k}\Phi_{\alpha k}|_{Y^{\tau}_{\alpha k} = (s/n)^{1/2}X + (1-s)^{1/2}\widehat{X}_{\alpha k}}\}|$$

$$\leq C_4(t, x)n^{-5/2}\left((\tau\sqrt{n})^5 + \mathbf{E}\left\{\left(\sum_{j \neq k}(X^{\tau}_{\alpha j})^2\right)^{5/2} + n^{3/2}\sum_j |\widehat{X}_{\alpha j}|^5\right\}\right)$$

$$\leq C_4(t, x),$$

where we took into account that we have, by the Hölder inequality and condition $|X_{\alpha j}^\tau| \le \tau \sqrt{n}$,

$$
\begin{aligned}
n^{-5/2} &\mathbf{E}\Big\{ \Big( \sum_{j \ne k} (X_{\alpha j}^\tau)^2 \Big)^{5/2} \Big\} \\
&\le n^{-5/2} \mathbf{E}^{5/6}\Big\{ \Big( \sum_j (X_{\alpha j}^\tau)^2 \Big)^3 \Big\} \\
&\le n^{-5/2} \Big( \sum_j \mathbf{E}\{ (X_{\alpha j}^\tau)^6 \} + 3 \sum_j \mu_{4,\alpha j}^\tau \sum_j \mu_{2,\alpha j}^\tau + \Big( \sum_j \mu_{2,\alpha j}^\tau \Big)^3 \Big)^{5/6} \\
&\le n^{-5/2} (\tau^2 n^2 \mu_4 + 3 n^2 \mu_4 a^2 + n^3 a^6)^{5/6} \le C < \infty
\end{aligned}
\tag{4.57}
$$

with $n$-independent $C$. We conclude that

$$
|\varepsilon_{3\tau}| \le C_4(t, x) \tau
\tag{4.58}
$$

[cf. (3.65)]. Besides, (4.56) and an analog of (3.30) for $X_{\alpha j}$ allow us to obtain for $T_{l\tau}$ of (4.49) an analog of (3.66) and (3.67):

$$
T_{l\tau} = T_l + r_l,
\tag{4.59}
$$

where now

$$
T_l = \frac{s^{(l-1)/2}}{l! n^{(l+1)/2}} \sum_{\alpha, k} \kappa_{l+1,\alpha k} \mathbf{E}\{ D_{\alpha k}^l(s) \Phi_{\alpha k} \},
\tag{4.60}
$$

$$
|r_l| \le C_l(t, x) \tau^{l-3} L_{mn}^{(4)}(\tau).
\tag{4.61}
$$

We have $T_0 = T_3 = 0$ (recall that $\kappa_{1,\alpha k} = \kappa_{4,\alpha k} = 0$), $T_1 = B_n$, and, in view of Lemma 4.1 below, $T_2 = o(1)$ [cf. (3.55)]. Hence,

$$
A_n = B_n + \varepsilon_{3\tau} + o(1),
$$

where the error term is a polynomial in $|t|$ and $|x|$ of degree 3 at most that vanishes as $m, n \to \infty$, $m/n \to c$ uniformly in $t$ and $x$ varying in a compact set $K \subset \{ t \ge 0, x \in \mathbb{R} \}$. This, (3.40), (4.47) and (4.58) imply (4.44) and complete the proof of the theorem. $\quad\square$

REMARK 4.2. A similar argument leads to the proof of the CLT for linear eigenvalue statistics of Hermitian analogs of (4.1)–(4.3), satisfying (1.5). The variance of the corresponding Gaussian law is $V_{\mathrm{Wish}}/2$, where $V_{\mathrm{Wish}}$ is given by (4.28). For real analytic test functions this formula is a particular case of the variance, obtained in [5] for random matrices $n^{-1} X^* T X$, where $X$ is a complex matrix with i.i.d. entries, satisfying (1.5) and (4.3), and $T$ is a certain Hermitian matrix.

LEMMA 4.1.   *We have under the conditions of Theorem 4.3*

$$(4.62) \qquad T_2 = \frac{s^{1/2}\mu_3}{2n^{3/2}} \sum_{\alpha,k} \mathbf{E}\{D_{\alpha k}^2(e_n^\circ(x)(Y^\tau U(t))_{\alpha k})\} = o(1)$$

*as* $m, n \to \infty$, $m/n \to c$.

PROOF.   By using (4.50)–(4.53), it can be shown that the assertion will follow from

$$T_{2p} = o(1), \qquad p = 1, 2, 3,$$

with

$$T_{21} = n^{-2} \sum_{\alpha,k} \mathbf{E}\{(X^\tau U)_{\alpha k}(t_1)\},$$

$$T_{22} = n^{-3} \sum_{\alpha,k} \mathbf{E}\{(X^\tau U)_{\alpha k}(t_1)(X^\tau U)_{\alpha k}(t_2)(X^\tau U)_{\alpha k}(t_3)\},$$

$$T_{23} = n^{-3} \sum_{\alpha,k} \mathbf{E}\{(X^\tau U X^{\tau T})_{\alpha\alpha}(t_1)(X^\tau U)_{\alpha k}(t_2)\}$$

[cf. (3.54)]. The Schwarz inequality, (2.16) and (4.3) yield

$$|T_{21}| = n^{-2}\left|\mathbf{E}\left\{\sum_j \left(\sum_k U_{jk}(t_1)\right)\left(\sum_\alpha X_{\alpha j}^\tau\right)\right\}\right|$$

$$\leq n^{-2}\mathbf{E}^{1/2}\left\{\sum_j \left|\sum_k U_{jk}(t_1)\right|^2\right\}\mathbf{E}^{1/2}\left\{\sum_j \left(\sum_\alpha X_{\alpha j}^\tau\right)^2\right\} \leq \sqrt{c_n} a n^{-1/2}$$

and [see also (4.57)]

$$|T_{22}| \leq n^{-5/2}\mathbf{E}^{1/2}\left\{\sum_{j_1,j_2,j_3}\left(\sum_\alpha X_{\alpha j_1}^\tau X_{\alpha j_2}^\tau X_{\alpha j_3}^\tau\right)^2\right\}$$

$$\leq n^{-5/2}\mathbf{E}^{1/2}\left\{\sum_j \left(\sum_\alpha (X_{\alpha j}^\tau)^3\right)^2 + 3\sum_{j_1 \neq j_2}\left(\sum_\alpha (X_{\alpha j_1}^\tau)^2 X_{\alpha j_2}^\tau\right)^2\right.$$

$$\left. + \sum_{j_1 \neq j_2 \neq j_3}\sum_\alpha (X_{\alpha j_1}^\tau)^2 (X_{\alpha j_2}^\tau)^2 (X_{\alpha j_3}^\tau)^2\right\}$$

$$\leq C(\tau^2 + 1)n^{-1/2}.$$

At last, by the Cauchy–Schwarz inequality, (2.16), (4.3) and (4.54), we have

$$|T_{23}| \leq n^{-3}\mathbf{E}^{1/2}\left\{\sum_\alpha |(X^\tau U X^{\tau T})_{\alpha\alpha}(t_1)|^2\right\}\mathbf{E}^{1/2}\left\{\sum_\alpha \left|\sum_k (X^\tau U)_{\alpha k}(t_2)\right|^2\right\}$$

$$\leq n^{-5/2}\mathbf{E}^{1/2}\left\{\sum_\alpha \left(\sum_j (X_{\alpha j}^\tau)^2\right)^2\right\}\mathbf{E}^{1/4}\left\{\sum_{j,j_1}\left(\sum_\alpha X_{\alpha j}^\tau X_{\alpha j_1}^\tau\right)^2\right\} \leq Cn^{-1/4}.$$

This completes the proof of the lemma.  □

4.4. *Central limit theorem for linear eigenvalue statistics of sample covariance matrices in the general case.*  Here we prove the CLT for the linear eigenvalue statistics of the sample covariance matrix not assuming that the fourth cumulant of its entries is zero [see (4.42)]. We use the scheme of the proof of Theorem 3.6 based on general differentiation formula (3.6) and the "a priory" bound (3.70) for the variance of statistics. Here is an analog of the bound for sample covariance matrices.

THEOREM 4.4.  *Let $M$ be the sample covariance matrix (4.1)–(4.3) satisfying (4.40) and (4.41), $M^\tau$ be corresponding truncated matrix (4.45), and*

$$u_{n\tau} = \operatorname{Tr} \exp\{itM^\tau\}.$$

*Then for any $\tau > 0$,*

$$(4.63) \qquad\qquad \mathbf{Var}\{u_{n\tau}(t)\} \le C_\tau(\mu_4)(1 + |t|^4)^2$$

*and*

$$(4.64) \qquad \mathbf{Var}\{\mathcal{N}_{n\tau}[\varphi]\} \le C_\tau(\mu_4)\left(\int (1 + |t|^4)|\widehat{\varphi}(t)|\, dt\right)^2,$$

*where $C_\tau(\mu_4)$ depends only on $\tau$ and $\mu_4$.*

We omit the proof of Theorem 4.4, because it repeats with natural modifications the proof of Theorem 3.5 for the Wigner case, and is again based on the use of the interpolation matrix (4.46) and known bound (4.29) for the Wishart matrix.

THEOREM 4.5.  *Let $M$ be the sample covariance matrix (4.1)–(4.3) satisfying (4.40) and (4.41), and $\varphi : \mathbb{R} \to \mathbb{R}$ be the test function satisfying (3.40). Then the centered linear eigenvalue statistic $\mathcal{N}_n^\circ[\varphi]$ of $M$ converges in distribution, as $m, n \to \infty$, $m/n \to c \in [1, \infty)$, to the Gaussian random variable of zero mean and variance*

$$(4.65) \quad V_{\mathrm{SC}}[\varphi] = V_{\mathrm{Wish}}[\varphi] + \frac{\kappa_4}{4c\pi^2 a^8}\left(\int_{a_-}^{a_+} \varphi(\mu)\frac{\mu - a_m}{\sqrt{4a^4 c - (\mu - a_m)^2}}\, d\mu\right)^2,$$

*where $V_{\mathrm{Wish}}[\varphi]$ is given by (4.28), $\kappa_4 = \mu_4 - 3a^4$ is the fourth cumulant of entries of $X$.*

PROOF.  Using the notation of the proof of Theorem 3.6, we note first that, according to Theorem 4.4 analogs of estimates (3.94) and (3.95), yielding the uniform boundedness of $Y_{n\tau}$ and $\partial Y_{n\tau}/\partial x$ remain valid in this case.

To estimate $\partial Y_{n\tau}/\partial t$, we write [cf. (3.96)]

$$\frac{\partial}{\partial t} Y_{n\tau}(x,t) = \frac{i}{\sqrt{n}} \sum_{\alpha,k} \mathbf{E}\{X_{\alpha k}^{\tau} \Phi_{\alpha k}\}$$

with

(4.66) $\quad \Phi_{\alpha k} = (Y^{\tau} U^{\tau})_{\alpha k}(t) e_{n\tau}^{\circ}(x), \qquad |\mathbf{E}\{D_{\alpha k}^{l} \Phi_{\alpha k}\}| \leq C_l(t,x), \qquad l \leq 5,$

and by using (3.6) and (4.59)–(4.61), we obtain an analog of (3.98):

$$\frac{\partial}{\partial t} Y_{n\tau}(x,t) = T_1 + O(1), \qquad T_1 = \frac{ia^2}{n} \sum_{\alpha,k} \mathbf{E}\{D_{\alpha k} \Phi_{\alpha k}\},$$

where, in view of (3.7), (4.61), (4.62) and (4.66), the error term is bounded by $C_2(t,x)$ in the limit (4.7). The term $T_1$ was calculated while deriving (4.33):

$$T_1 = ia^2 c_n Y_{n\tau}(x,t) + ia^2 t \mathbf{E}\{n^{-1} u_{n\tau}'(t) e_{n\tau}^{\circ}(x)\}$$

$$+ ia^2 \int_0^t \mathbf{E}\{n^{-1} u_{n\tau}'(t-t_1) u_{n\tau}(t_1) e_{n\tau}^{\circ}(x)\} \, dt_1$$

$$- 2a^2 x \int t_1 \widehat{\varphi}(t_1) \mathbf{E}\{n^{-1} u_{n\tau}'(t+t_1) e_{n\tau}(x)\} \, dt_1.$$

We also have, by (2.16) and (4.3),

(4.67)
$$\mathbf{E}\{|n^{-1} u_{n\tau}'(t)|^2\} = n^{-2} \mathbf{E}\{|\operatorname{Tr} M^{\tau} U^{\tau}(t)|^2\}$$

$$\leq n^{-3} \mathbf{E}\left\{ \sum_{j,k} \left( \sum_{\alpha} X_{\alpha j}^{\tau} X_{\alpha k}^{\tau} \right)^2 \right\} \leq C$$

and, by integrating by parts,

$$\int_0^t \mathbf{E}\{n^{-1} u_{n\tau}'(t-t_1) u_{n\tau}(t_1) e_{n\tau}^{\circ}(x)\} \, dt_1$$

$$= \int_0^t \mathbf{E}\{n^{-1} u_{n\tau}'(t-t_1)\} Y_{n\tau}(x,t_1) \, dt_1$$

$$+ \int_0^t \mathbf{E}\{u_{n\tau}^{\circ}(t-t_1) n^{-1} u_{n\tau}'(t_1) e_{n\tau}^{\circ}(x)\} \, dt_1,$$

where the r.h.s. is uniformly bounded in view of (4.63) and (4.67). Hence, $T_1$ is uniformly bounded for any $\tau > 0$, and so does $\partial/\partial t Y_{n\tau}$. This and analogs of (3.94) and (3.95) imply the existence of a subsequence $\{Y_{n_l \tau}\}$ that converges uniformly to a continuous $Y_{\tau}$.

Now an argument similar to that leading to (3.104)–(3.106) yields an analog of (3.104):

$$Y_{n\tau}(x,t) = T^\tau_{a^2,n} + T^\tau_{\kappa_4,n} + \mathcal{E}_{3\tau,n}(t,x) + o(1),$$

(4.68)

$$n, m \to \infty, m/n \to c,$$

where the terms on the r.h.s. are given by the r.h.s. of (3.105)–(3.107) with $a^2$ instead of $w^2(1+\delta_{jk})$, $\Phi_{\alpha k}$ of (4.66), and for any $\tau > 0$ the reminder term $o(1)$ vanishes in the limit (4.7) uniformly on any compact set of $\{t \geq 0, x \in \mathbb{R}\}$.

The term $T^\tau_{a^2,n}$ was in fact calculated in the proof of Theorem 4.2 and is equal to $Y_n(x,t)$ of (4.34) and (4.35) with the Wishart matrix $M$ replaced by the truncated sample covariance matrix $M^\tau$. Using (4.63) to estimate the reminder term $r_n$ of (4.35), and noting that by an analog of (3.58) $v_{n\tau} \to v_{MP}$ in the limit (4.7), we get an analog of (3.108) in the same limit:

$$T^\tau_{a^2,n} \to ia^2(c-1) \int_0^t Y_\tau(x,t_1)\, dt_1$$

$$+ 2ia^2 \int_0^t v_{MP}(t-t_1) Y_\tau(x,t_1)\, dt_1 + x Z_\tau(x) A(t)$$

with $A(t)$ defined in (4.37).

Consider now the term $T^\tau_{\kappa_4,n}$ of (4.68), given by (3.106) with $\Phi_{\alpha k}$ of (4.66). It follows from (4.50)–(4.53) and an argument similar to that of the proof of Lemma 4.1 that the contribution to $T^\tau_{\kappa_4,n}$ due to any term of $n^{-2}\sum_{\alpha,k} D^3_{\alpha k}\Phi_{\alpha k}$, containing at least one element $(Y^\tau U^\tau)_{\alpha k}$, vanishes as $m, n \to \infty$, $m/n \to c$. Thus, we are left with the terms, containing only diagonal elements of $U^\tau$ and $Y^\tau U^\tau Y^{\tau T}$. These terms arise from $e^\circ_{n\tau} D^3_{\alpha k}(Y^\tau U^\tau)_{\alpha k}$ and $3D_{\alpha k}(Y^\tau U^\tau)_{\alpha k} D^2_{\alpha k} e^\circ_{n\tau}$, and by (4.50)–(4.53), their contributions to $T^\tau_{\kappa_4,n}$ are [cf. (3.110) and (3.111)]

$$-\frac{\kappa_4}{n^2} \int_0^t \sum_{\alpha,k} \mathbf{E}\{([U^\tau_{kk} + i(Y^\tau U^\tau Y^{\tau T})_{\alpha\alpha} * U^\tau_{kk}]$$

(4.69)

$$* [U^\tau_{kk} + i(Y^\tau U^\tau Y^{\tau T})_{\alpha\alpha} * U^\tau_{kk}])(t_1) e^\circ_{n\tau}(x)\}\, dt_1$$

and

$$-\frac{ix\kappa_4}{n^2} \sum_{\alpha,k} \mathbf{E}\Big\{ e_{n\tau}(x) \int_0^t (U^\tau_{kk} + i(Y^\tau U^\tau Y^{\tau T})_{\alpha\alpha} * U^\tau_{kk})(t_1)\, dt_1$$

(4.70)

$$\times \int t_2 \widehat{\varphi}(t_2)(U^\tau_{kk} + i(Y^\tau U^\tau Y^{\tau T})_{\alpha\alpha} * U^\tau_{kk})(t_2)\, dt_2 \Big\}.$$

Thus, the entries of $U^\tau$ and $Y^\tau U^\tau Y^{\tau T}$ are present here in the form [cf. (3.112) and (3.113)]

(4.71)          $$K_{p0} = \mathbf{E}\{v_{n\tau}(t_1,t_2) w_{p,n}(t_3,t_4) e^\circ_{n\tau}(x)\}, \qquad p = 0, 1,$$

and

(4.72) $\qquad K_p = \mathbf{E}\{v_{n\tau}(t_1,t_2)w_{p,n}(t_3,t_4)e_{n\tau}(x)\}, \qquad p = 0,1,$

where $v_{n\tau}(t_1,t_2)$ is defined analogously to (3.114) and satisfies an analog of (3.115), and

(4.73) $\qquad w_{p,n}(t_3,t_4) = n^{-1}\sum_{\alpha}(YUY^{\tau T})_{\alpha\alpha}(t_3)(YUY^{\tau T})^p_{\alpha\alpha}(t_4)$

satisfies

$$|\mathbf{E}\{w_{p,n}(t_3,t_4)\}| \le n^{-(2+p)}\mathbf{E}\left\{\sum_{\alpha}\left(\sum_j (X^\tau_{\alpha j})^2\right)^{(1+p)}\right\} \le C$$

by (4.3) and (4.54). Since the expectations of $v_{n\tau}(t_1,t_2)$ and $w_{p,n}(t_3,t_4)$ are uniformly bounded, and by Lemma 4.2 below their variances vanish in subsequent limit (4.43), then, applying the Schwarz inequality and (3.116), we conclude that

(4.74) $\qquad \begin{aligned} &K_{p0} = o(1),\\ &K_p = Z_{n\tau}(x)\overline{v}_{n\tau}(t_1,t_2)\mathbf{E}\{w_{p,n}(t_3,t_4)\} + o(1), \qquad p = 0,1 \end{aligned}$

[cf. (3.118) and (3.119)], where the error terms vanish in the limit (4.7) uniformly in $(t,x) \in \mathbb{R}^2$. Using the interpolation argument similar to that in the proof of Lemma 3.1 with the GOE matrix replaced by the Wishart matrix, we get an analog of (3.120):

(4.75) $\qquad \lim_{m,n\to\infty,m/n\to c}\overline{v}_{n\tau}(t_1,t_2) = v_{MP}(t_1)v_{MP}(t_2).$

To find the limit of $\mathbf{E}\{w_{0,n}(t_3,t_4)\}$, we note that

$$\mathbf{E}\{w_{0,n}(t_3,t_4)\} = i^{-1}\overline{v}'_{n\tau}(t_3),$$

where $\overline{v}_{n\tau}(t)$ converges to $v_{MP}(t)$ as $m,n\to\infty$, $m/n\to c$, and that by (4.67) and a similar argument, the sequences $\{\overline{v}'_{n\tau}(t)\}$ and $\{\overline{v}''_{n\tau}(t)\}$ are uniformly bounded, so that we have

$$\lim_{m,n\to\infty,m/n\to c}\mathbf{E}\{\overline{v}'_{n\tau}(t)\} = i^{-1}v'_{MP}(t)$$

uniformly in $t$, varying in a finite interval. Furthermore, it can be shown by an argument, used not once before and based on (2.14), (3.6) and relation (4.80) below, that the functions

$$\lim_{m,n\to\infty}\mathbf{E}\{(Y^\tau U^\tau(Y^\tau)^T)_{\alpha\alpha}(t)\}$$

and

$$\lim_{m,n\to\infty}\mathbf{E}\left\{m^{-1}\sum_{\alpha}(Y^\tau U^\tau(Y^\tau)^T)_{\alpha\alpha}(t)\right\} = (ic)^{-1}v'_{MP}(t)$$

satisfy the integral equation

$$h(t) = a^2 v_{MP}(t) + ia^2 \int_0^t h(t - t_1) v_{MP}(t_1) \, dt_1.$$

This and Proposition 2.1(v) imply that the functions coincide, and we obtain, in view of (4.80),

$$\lim_{m,n \to \infty, m/n \to c} \mathbf{E}\{w_{1,n}(t_3, t_4)\} = -c^{-1} v'_{MP}(t_3) v'_{MP}(t_4).$$

We conclude from the above that the contribution of (4.69) to $T_{\kappa_4, n_l}$ vanishes as $n_l \to \infty$ uniformly in $t$ and $x$, varying in any compact set of $\{t \geq 0, x \in \mathbb{R}\}$, while in (4.70) we can replace $U_{kk}$ by $v_{MP}$ and $(Y^\tau U^\tau Y^{\tau T})_{\alpha\alpha}$ by $(ic)^{-1} v'_{MP}$. We obtain

$$(4.76) \qquad \lim_{m,n \to \infty, m/n \to c} T^\tau_{\kappa_4, n} = -c^{-1} \kappa_4 x Z_\tau(x) C[\varphi] \int_0^t A_{\kappa_4}(t_1) \, dt_1,$$

where

$$(4.77) \qquad A_{\kappa_4}(t) = c v_{MP}(t) + \int_0^t v_{MP}(t - t_1) v'_{MP}(t_1) \, dt_1,$$

$$(4.78) \qquad C[\varphi] = i \int t \widehat{\varphi}(t) A_{\kappa_4}(t) \, dt$$

or, in view of Proposition 2.1, (4.20), (4.21) and (4.38) [cf. (2.81)],

$$A_{\kappa_4}(t) = \frac{1}{2\pi a^4} \int_{a_-}^{a_+} e^{i\mu t} \sqrt{4a^4 c - (\mu - a_m)^2} \, d\mu.$$

Plugging the last expression in (4.78) and integrating by parts, we get

$$C[\varphi] = \frac{1}{2\pi a^4} \int_{a_-}^{a_+} \varphi(\mu) \frac{\mu - a_m}{\sqrt{4a^4 c - (\mu - a_m)^2}} \, d\mu.$$

This, (4.68) and (4.76) lead to the integral equation for $Y_\tau(x, t)$ [cf. (3.127) and (3.128)]:

$$Y_\tau(x, t) - ia^2(c - 1) \int_0^t Y_\tau(x, t_1) \, dt_1 - 2ia^2 \int_0^t v_{MP}(t - t_1) Y_\tau(x, t_1) \, dt_1$$

$$= -x Z_\tau(x) \left( A(t) + \kappa_4 c^{-1} C[\varphi] \int_0^t A_{\kappa_4}(t_1) \, dt_1 \right) + \mathcal{E}_{3\tau}(t, x),$$

where $\mathcal{E}_{3\tau}$ satisfies (3.127).

Now, to finish the proof, we have to follow the part of the proof of Theorem 3.6 after (3.127) to obtain (4.65).  $\square$

LEMMA 4.2. *We have, under the conditions of Theorem 4.5 in the limit (4.7),*

$$\textbf{Var}\{U^\tau_{kk}(t)\} = o(1), \tag{4.79}$$

$$\textbf{Var}\{(Y^\tau U^\tau Y^{\tau T})_{\alpha\alpha}(t)\} = o(1). \tag{4.80}$$

PROOF. The proof of (4.79) repeats with natural modifications the one of an analogous assertion for the Wigner matrix (see Lemma 3.1). It is based on the interpolation procedure and following from the Poincaré inequality (2.21) validity of (4.79) for Wishart matrices.

To prove (4.80), we consider

$$V_{n\tau}(t_1, t_2) = \textbf{E}\{(Y^\tau U^\tau (Y^\tau)^T)_{\alpha\alpha}(t_1)(Y^\tau U^\tau (Y^\tau)^T)^\circ_{\alpha\alpha}(t_2)\},$$

putting in an appropriate moment $t_2 = -t_1$ to get $\textbf{Var}\{(Y^\tau U^\tau (Y^\tau)^T)_{\alpha\alpha}(t_1)\}$. We have, by (3.6) and (4.56),

$$V_{n\tau}(t_1, t_2) = n^{-1/2} \sum_k \textbf{E}\{X_{\alpha k}(Y^\tau U^\tau)_{\alpha k}(t_1)(Y^\tau U^\tau (Y^\tau)^T)^\circ_{\alpha\alpha}(t_2)\}$$

$$= n^{-1} \sum_k \textbf{E}\{D_{\alpha k}((Y^\tau U^\tau)_{\alpha k}(t_1)(Y^\tau U^\tau (Y^\tau)^T)^\circ_{\alpha\alpha}(t_2))\} + \varepsilon_{1\tau},$$

where $|\varepsilon_{1\tau}| \le \mu_4 C_2(t, x)$ by an argument similar to that used in (4.57) and, by (4.52) and (4.53), the sum on the r.h.s. is

$$i \int_0^{t_1} \overline{v}_{n\tau}(t_1 - s) V_{n\tau}(s, t_2)\, ds + \textbf{E}\{v_{n\tau}(t_1)(Y^\tau U^\tau (Y^\tau)^T)^\circ_{\alpha\alpha}(t_2)\}$$

$$+ i \int_0^{t_1} \textbf{E}\{v^\circ_{n\tau}(t_1 - s)(Y^\tau U^\tau Y^{\tau T})_{\alpha\alpha}(s)(Y^\tau U^\tau (Y^\tau)^T)^\circ_{\alpha\alpha}(t_2)\}\, ds$$

$$+ it_1 n^{-1} V_{n\tau}(t_1, t_2) + 2n^{-1}\textbf{E}\{(Y^\tau U^\tau Y^{\tau T})_{\alpha\alpha}(t_1 + t_2)\}$$

$$+ 2in^{-1} \int_0^{t_2} V_{n\tau}(t_1 + s, t_2 - s)\, ds.$$

It follows from (4.3) and (4.54) that $V_{n\tau}(t_1, t_2)$ and $\textbf{E}\{(Y^\tau U^\tau (Y^\tau)^T)^2_{\alpha\alpha}(t)\}$ are uniformly bounded. This, the Schwarz inequality and (4.63) imply that

$$\textbf{E}\{v_{n\tau}(t_1)(Y^\tau U^\tau (Y^\tau)^T)^\circ_{\alpha\alpha}(t_2)\} = O(n^{-1})$$

and

$$|\textbf{E}\{v^\circ_{n\tau}(t_1 - s)(Y^\tau U^\tau (Y^\tau)^T)_{\alpha\alpha}(s)(Y^\tau U^\tau (Y^\tau)^T)^\circ_{\alpha\alpha}(t_2)\}|$$

$$\le \textbf{E}\left\{|v^\circ_{n\tau}(t_1 - s)|\left(\sum_j (Y^\tau_{\alpha j})^2\right)^2\right\}$$

$$\leq \mathbf{E}^{1/4}\{|v_{n\tau}^\circ(t_1-s)|^2\}\mathbf{E}^{1/4}\Big\{|v_{n\tau}^\circ(t_1-s)|^2\Big(\sum_j(Y_{\alpha j}^\tau)^2\Big)^2\Big\}$$

$$\times\mathbf{E}^{1/2}\Big\{\Big(\sum_j(Y_{\alpha j}^\tau)^2\Big)^3\Big\}=O(n^{-1/2})$$

as $m,n\to\infty$, $m/n\to c$. Besides, a bit tedious but routine calculations, similar to those in the proof of Lemma 4.1, yield the boundedness of derivatives of $V_{n\tau}(t_1,t_2)$ for any $\tau>0$. Thus, there exists a subsequence $(m_l,n_l)$ such that the limit $V(t_1,t_2)=\lim_{m_l,n_l\to\infty}V_{n_l,\tau}(t_1,t_2)$ exists and satisfies the equation

$$V(t_1,t_2)=i\int_0^{t_1}\overline{v}_{MP}(t_1-s)V(s,t_2)\,ds.$$

Now Proposition 2.1 implies that $V(t_1,t_2)=0$. This completes the proof of the lemma. □

Mathematical Division
B. Verkin Institute
  for Low Temperature Physics
47 Lenin Avenue
61103, Kharkov
Ukraine
E-mail: lytova@ilt.kharkov.ua
        lpastur@flint.ilt.kharkov.ua